\documentclass[10pt]{amsart}

\usepackage{amssymb,amsmath,amscd,amsfonts,a4,psfrag,graphicx,
latexsym,epsfig,color}
\usepackage[active]{srcltx}
\usepackage[all]{xy}
\usepackage{tikz}
\usetikzlibrary{matrix} \everymath{\displaystyle}

\newtheorem{teo}{Theorem}[section]
\newtheorem{lem}[teo]{Lemma}
\newtheorem{cor}[teo]{Corollary}

\newtheorem{prop}[teo]{Proposition}
\newtheorem{defi}[teo]{Definition}
\newtheorem{ques}[teo]{Question}
\newtheorem{conj}[teo]{Conjecture}
\newtheorem{remark}[teo]{Remark}

\newcommand{\cvd}{\hfill$\Box$}

\newcommand{\R}{\mathbb{R}}

\newcommand{\Z}{\mathbb{Z}}

\newcommand{\Aa}{{\mathcal A}}
\newcommand{\Bb}{{\mathcal B}}

\newcommand{\Ff}{{\mathcal F}}

\newcommand{\Hh}{{\mathcal H}}

\newcommand{\Mm}{{\mathcal M}}
\newcommand{\Nn}{{\mathcal N}}

\newcommand{\Pp}{{\mathcal P}}

\newcommand{\Ss}{{\mathcal S}}
\newcommand{\Tt}{{\mathcal T}}
\newcommand{\Uu}{{\mathcal U}}
\newcommand{\Ww}{{\mathcal W}}
\newcommand{\Vv}{{\mathcal V}}

\newcommand{\PP}{{\bf P}}

\newcommand{\bG}{{\mathfrak b}}

\newcommand{\SG}{{\mathfrak S}}

\newcommand{\fG}{{\mathfrak f}}
\newcommand{\mG}{{\mathfrak m}}

\newcommand{\Dim}{{\it Proof.\ }}

\def\cvd{\hfill$\Box$}

\title{Ideal triangulations of $3$-manifolds up to decorated transit equivalences }
\begin{document}
\author{Riccardo Benedetti} \address{Dipartimento di Matematica, Largo
  Bruno Pontecorvo 5, 56127 Pisa, Italy}
\email{riccardo.benedetti@unipi.it} 

\begin{abstract}
We consider $3$-dimensional pseudo-manifolds $\hat M$ with a given set
of marked point $V$ such that $\hat M \setminus V$ is the interior of
a compact $3$-manifold with boundary. An ideal triangulation $T$ of
$(\hat M, V)$ has $V$ as set of vertices. A branching $(T,b)$ enhances
$T$ to a $\Delta$-complex. Branched triangulations of $(\hat M, V)$
are considered up to the $b$-transit equivalence generated by isotopy
and ideal branched moves which keep $V$ pointwise fixed. We extend a
well known connectivity result for `naked' ideal triangulations by
showing that branched ideal triangulations of $(\hat M,V)$ are
equivalent to each other.  A pre-branching $(T,\omega)$ is a system of
transverse orientations at the $2$-facets of $T$ verifying a certain
global constraint; pre-branchings are considered up to a natural
$pb$-transit equivalence.  If $M$ is oriented, every branching $(T,b)$ induces a
pre-branching $(T,\omega_b)$ and every $b$-transit induces a
$pb$-transit.  The quotient set of pre-branchings up to transit
equivalence is far to be trivial; we get some information about it and
we characterize the pre-branchings of the type $\omega_b$.
Pre-branched and branched moves are naturally organized in subfamilies
which give rise to restricted transit equivalences.  In the branching
setting we revisit, with some complement, early results about the {\it
  sliding} transit equivalence and outline a (partially conjectural)
conceptually different approach to the branched ideal connectivity and
eventually also to the naked one. The basic idea is to point out some
structures of differential topological nature on $M$ which are carried
by every branched ideal triangulation $(T,b)$ of $\hat M$, are
preserved by the sliding transits and can be modified by the full
branched transits.  The {\it non ambiguos} transit equivalence already
widely studied on pre-branchings lifts to a specialization of the
sliding equivalence on branched triangulations; we point out a few
specific insights, again in terms of carried structures preserved by
the non ambiguous and which can be modified by the whole sliding
transits.

\end{abstract}

\maketitle

\section{Introduction}
This paper concerns $3$-manifold triangulations.  The recent article
\cite{RST} is a fresh and valuable reference for many results about
``naked'' triangulations (with a lot of expressive pictures).  We
widely refer to its body and bibliography.  We work on a given compact
connected smooth $3$-manifold $M$ with {\it non empty} boundary
$\partial M$.  We denote by $\hat M$ the space obtained by collapsing
to one point $v$ each boundary component, we denote $V$ the set of
these points. Sometimes $\hat M$ is said a {\it pseudomanifold}.  Then
the interior Int$(M)$ is embedded in $\hat M$, onto $\hat M\setminus
V$.  The non manifold points of $\hat M$ are the points of $V$
associated to non spherical components of $\partial M$.  We consider
possibly loose triangulations of $\hat M$ such that the set of
vertices coincides with $V$. Sometimes a vertex which is a manifold
point is called a {\it material vertex}.  ``Loose'' means that self
and multiple face adjiacency are allowed. Such a triangulation is
usually called an {\it ideal triangulation} of Int$(M)$, understanding
that the vertices are ``at infinity''.  This terminology alludes to
the triangulations by ideal tetrahedra of a hyperbolic $3$-manifold
with cusps. However we simply call them ideal triangulations of $\hat
M$.  It is sometimes useful to consider an ideal triangulation as a
way to realize $\hat M$ by assembling ``abstract'' tetrahedra by
gluing their abstract $2$-faces in pairs in such a way that no face
remains unglued. Every $\hat M$ admits ideal triangulations.  The
ideal triangulations of $\hat M$ are considered up to the {\it ideal
  transit equivalence} which is generated by two basic local moves
(Section \ref{generality}) and isotopy which keep $V$ pointwise fixed.
These basic {\it ideal moves}
are  the $2 \leftrightarrow 3$ 
and the {\it quadrilateral} $0
\leftrightarrow 2$ {\it move}.
The numbers  in the name refers
to the variation of the number of tetrahedra when the move is performed; the move
is {\it positive} if this number increases. Denote by
$\Tt^{id}(M)$ the corresponding quotient set.

The {\it completed transit equivalence} is obtained by adding one
further move called {\it triangular} $0 \leftrightarrow 2$ {\it
  move}. Equivalently we can add instead the {\it stellar} $1
\leftrightarrow 4$ {\it move}. In fact we will use freely both
moves. After a positive such a move we no longer have an ideal
triangulation of $\hat M$. Instead we have an ideal triangulation of
$\hat M'$, where $M'$ is obtained from $M$ by removing the interior of
a three ball embedded into Int$(M)$. So there is a new material vertex
in $\hat M'$.  If such a positive move occurs in a composite transit
$T_1\Rightarrow T_2$ (relative isotopy will be always understood)
connecting two ideal triangulations of $\hat M$, then it must be
compensated later by a negative inverse move. We denote by $\Tt(M)$
the quotient set under the completed equivalence.

{\bf The dual viewpoint.}  For every ideal triangulation $T$ of $\hat
M$, the $2$-skeleton $\Sigma=\Sigma_T$ of the dual cell decomposition
is a {\it standard} (internal) spine of $M$. Here ``standard'' means
that the spine has generic singularities and every stratum of its
natural stratification is an open cell of the appropriate dimension.
Sometimes one says ``special'' instead of ``standard''. If we drop out
the cellularity condidion, then we have the notion of ``simple''
spine.  A local portion of $\Sigma$ corresponding to a tetrahedron is
called a {\it butterfly}; the fully symmetric picture of it is the
cone based on the $1$-skeleton of a tetrahedron with centre at an
interior point - one sees it on the left side of Figure \ref
{butterfly}.  The $1$-skeleton of $\Sigma$ is its singular set
Sing$(\Sigma)$. Every $0$-cell of Sing$(\Sigma)$ (also called a
vertex) is quadrivalent.  The above moves on ideal triangulations can
be fully rephrased in terms of standard spines.  The spine version of
the quadrilateral $0\leftrightarrow 2$ move is called {\it lune}
move. The spine version of the triangular $0\leftrightarrow 2$ move is
called {\it bubble} move. We will freely adopt both equivalent dual
viewpoints.

\begin{remark}{\rm Although they are equivalent, there is some
    qualitative difference beetween standard spines and ideal
    triangulations. For example a triangulation $2\to 3$ move
    basically is a {\it discrete} transition with a cell decomposition
    as intermediate ``state'' which is no longer a triangulation. The
    corresponding spine transition can be realized by a {\it
      continuous} deformation passing through a non generic spine. By
    this difference, sometimes manipulation of spines is easier as
    more visual. }
\end{remark}

\smallskip

{\bf The naked connectivity results.}  The following are fundamental
well known connectivity results for ``naked'' triangulations.

\begin{teo}\label{nakedC} The set $\Tt(M)$ consists of one point.
\end{teo}  

\begin{teo}\label{nakedI} The set $\Tt^{id}(M)$ consists of one point.
\end{teo}
\smallskip

A few comments are in order:
\begin{enumerate}
\item These theorems are  definitely {\it weaker} versions of the connectivity results
  discussed in \cite{RST} (Theorem 1.2 and Theorem 1.4 therein
  respectively).  The {\it strong} versions are obtained by discarding the
  quadrilateral $0\leftrightarrow 2$ move from the generators of the
  completed relation, and also from the ideal relation, provided that
  one deals with triangulations with at least two tetrahedra.

\item In both weak and strong versions, Theorem \ref{nakedC} obviously
  is an easy consequence of Theorem \ref{nakedI}. However, the
  current proof
   discussed in \cite{RST} of the strong version of the second is based on the
  validity of the strong version of the first.

\item In \cite{RST} one can find an accurate discussion about the
  contributions of several authors to the proof of the strong versions
  of the theorems. The
  essential ideas to derive the second  from the first are due
  to Matveev \cite{ Mat} and Piergallini \cite{P}; a detailed implementation
  as a particular case of a more general result is in \cite{Am}.

\item
  The present weak version of the theorems is adequate to our aims for
  several reasons:

  (a) Since our early motivations by dealing with the quantum
  hyperbolic invariants, in the decorated setting we are mainly
  concerned with the opposition {\it ideal vs completed} transit
  equivalence, not with the search of a minimal set of generating
  moves; for example the quadrilateral $0\leftrightarrow 2$ move
  occurs very naturally in the treatment of non ambiguous structures
  carried by layered triangulations of $3$-manifolds fibred on $S^1$
  with punctured fibre and their reduced QHI, see \cite{NA},
  \cite{QT}, \cite{A}.
  
(b) More substantially, we assume a refinement of Theorem \ref{nakedI}
  (hence of \ref{nakedC}) in our weak version, due to \cite{Ma}; for
  every ideal triangulations $T_1$ and $T_2$ of $\hat M$, this
  provides the existence of composite ideal transits towards a {\it
    same} triangulation $T$, $T_1\Rightarrow T$, $T_2\Rightarrow T$,
  both enterely composed by {\it positive} ideal moves ( $2\rightarrow
  3$ {\it and} $0\rightarrow 2$). Positive moves behave well with
  respect to their decorated enhancements; see in particular Lemma
  \ref{sameT} which is our starting point to treat decorated transit
  equivalences.
\end{enumerate} 

\smallskip

{\bf Decorated triangulations and their transit equivalences.}
Different notions of decorated ideal triangulations of $3$-manifolds
considered up to suitable transit equivalences naturally arise for
example in the developments of {\it quantum hyperbolic geometry}, see
in particular \cite{AGT}, \cite{NA}, \cite{QT}, and in several other
instances of quantum invariants based on some ``$6j$ symbols'' theory
(see \cite{NA}); on the other hand, the theme of the combinatorial
realization of different kind of structures on $M$ of differential
topological nature in terms of suitably decorated triangulations up to
transit has been faced since \cite{BP2} to \cite{BP}, see also
\cite{AFJ}, \cite{FJR}, \cite{MM}.  Here we consider two basic kinds
of decoration.  Let $T$ be an ideal triangulation of $\hat M$ and
$\Sigma$ be the dual spine of $M$.
\begin{itemize}
\item {\bf Branching:} A branching $b$ on $T$ is a system of
  orientations of the edges of $T$ which lifted to every abstract
  tetrahedron $(\Delta,b)$ of $T$ is induced by a (local) ordering of
  the vertices, so that every edge goes towards the bigest
  endpoint. Equivalently $b$ enhances $T$ to be a $\Delta$-{\it
    complex} accordingly with \cite{HATCHER}, Chapter 2. Usually an
  order of the verices of $\Delta$ will be specified by a labelling by
  $0,1,2,3$. Dually $(T,b)$ becomes $(\Sigma, b)$ (we keep the
  notation $b$) that is a system of transverse orientations of the
  open $2$-cells of $\Sigma$ (also called regions). 
  \smallskip
  
  {\it $M$ oriented:}  In this case $(\Sigma,b)$ is equivalent to a 
  system of orientations of the
  regions, by stripulating that every $b$-oriented edge of $(T,b)$
  intersects its dual oriented region with intersection number equal
  to $1$.  Then the branching condition translates in the property
  that at every open $1$-cell $e$ in Sing$(\Sigma)$ (also called an
  edge) the three incident local branches of oriented regions of
  $(\Sigma,b)$ do not induce the same orientation on $e$; hence there
  is a {\it prevailing orientation} induced by two of them. 
  If $M$ is
  oriented, to every branched ideal triangulation $(T,b)$ is associated
  a simplicial fundamental $3$-cycle
  $$Z(T,b):= \sum_{\Delta \subset T} *_{(\Delta,b)} (\Delta,b)$$
  where $(\Delta,b)$ varies among the $3$-simplexes of the $\Delta$-complex
  $(T,b)$, the sign $*_{(\Delta,b)}=\pm 1 \in \Z$ and equals $1$ if and only if
  the $b$-orientation of $\Delta$ agrees with the ambient orientation of $M$.
  
  \begin{figure}[ht]
\begin{center}
 \includegraphics[width=7cm]{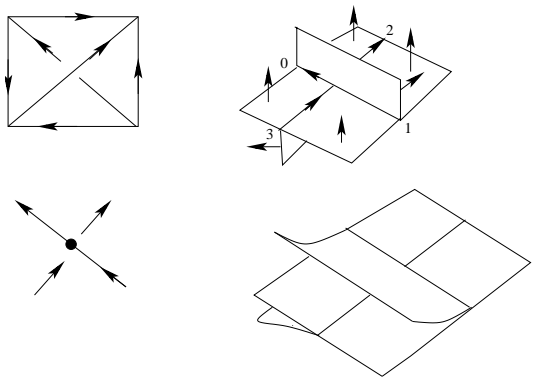}
\caption{\label{b+} From a branched tetrahedron to branched butterfly.}
\end{center}
\end{figure}
 
 In every case, the branching encodes a structure of smooth {\it
   transversely oriented branched surface} on $\Sigma$ embedded in the
 interior of $M$; hence $\Sigma$ becomes a so called {\it branched
   standard spine} of $M$ (see \cite{BP2} for all details). This also
 justifies the name ``branching''.  A branched butterfly associated to
 a branched tetrahedron is shown in Figure \ref{b+}; the labels
 $0,1,2,3$ refer to the order on the dual $2$-faces of $(\Delta,b)$
 determined by the branching (accordingly with the order of the
 opposite vertices); the decorated planar crossing on the left-bottom
 of the figure encodes all the other pictures.
  
\smallskip  
  
\item{\bf Pre-branching:}  A pre-branching $\omega$ on $T$ is a system
  of orientations of the open $1$-cells of Sing$(\Sigma)$ such that at
  every quadrivalent vertex of Sing$(\sigma)$ two branches are ingoing and two are
  outgoing. This can be equivalently rephrased in terms of a system of
  transverse orientations on the $2$-faces of $T$.
  \end{itemize}  

\smallskip

Given a branched triangulation $(T,b)$ or a pre-branched triangulation
$(T,\omega)$ of $\hat M$ every {\it positive} naked move $T\to T'$
(ideal or not) can be enhaced to some decorated move $(T,b)\to
(T',b')$ or $(T,\omega)\to (T',\omega')$. This means that in both
cases the decorations coincide on the portion of triangulation which
persists in both $T$ and $T'$; such a decorated move is also called a
$b$- or $pb$-{\it transit}.  A negative decorated move is by
definition the inverse of a positive one.  

If $M$ is {\it oriented}, every branching $b$ induces
a pre-branching $\omega_b$, by taking the prevailing dual edge orientation as above;
every $b$-move $b\to b'$ induces a
$pb$-move $\omega_b \to \omega_{b'}$.
\smallskip

By using only ideal moves we define the {\it decorated ideal transit
  equivalences} with quotient sets denoted $\Bb^{id}(M)$ and
$\Pp\Bb^{id}(M)$ respectively.  Similarly we define the {\it
  completed} decorated relations with quotient sets $\Bb(M)$ and
$\Pp\Bb(M)$.  Pre-branchings are much more flexible than branchings
and are the basic ingredient for the theory of {\it taut}
triangulations \cite{La}.  In particular every naked ideal
triangulation $T$ carries pre-branchings while there are naked
triangulations, even for oriented $M$, which does {\it not} carry any
branching.  On the other hand, every $M$ admits branched ideal
triangulations.  Quantum hyperbolic state sums were originally
supported by branched triangulations (of oriented manifolds); since
\cite{AGT} it is clear that pre-branching is definitely the most
funtamental structure governing such state sums. Often pre-branchings
(on oriented $M$) are encoded by auxiliary decorations called {\it
  weak-branchings} which reveal in a more trasparent way a geometric
content. We will occasionally make use of them referring for the
details to \cite{NA}, \cite{AGT} and also \cite{BP}.
 
 \smallskip

{\bf Branched connectivity results.} 
Here is two main results of the
paper, see Section \ref{b-connectivity}.

\begin{teo}\label{bC} For every $M$, $\Bb(M)$ consists of one point.
\end{teo}

\begin{teo}\label{bI} For every $M$, $\Bb^{id}(M)$  consists of one point.
\end{teo}

A few comments are in order:
\begin{enumerate}
\item A somewhat demanding proof of Theorem \ref{bC} is given in
  \cite{Cost}, starting from Lemma \ref{sameT}, in terms of branched spines.
  We will provide an
  alternative quick proof in terms of triangulations, based again on Lemma
  \ref{sameT} and the first two steps in the realization of the
  (naked) barycentric subdivision of a given ideal triangulation $T$
  by means of positive $2\to 3$ and $1\to 4$ moves (see \cite{RST}
  Section 2.5). A simple but important ingredient will be the
  elementary move of {\it inverting a good ambiguous edge} which can
  be added without modifying the $b$-ideal transit equivalence.

\item As for the naked ideal connectivity result in the strong version
  mentioned above, we will give a proof of Theorem \ref{bI} based on
  the validity of Theorem \ref{bC}, by applying suitable branched
  versions of Matveev's {\it arch and associated constructions}.  The
  underlying naked configurations which we actually apply these
  constructions to are considerably simpler that the general ones
  faced in the proof of the strong version of Theorem \ref{nakedI};
  also the branched enhancement will benefit from this
  simplification.
\end{enumerate}

\smallskip

{\bf On the structure of $\Pp\Bb^{id}$.} In contrast with
$\Bb^{id}(M)$, in general $\Pp\Bb^{id}(M)$ is far to be trivial, even
infinite. We will get some information about it, in particular we
characterize the (one-point) image of $\Bb^{id}(M)$ in
$\Pp\Bb^{id}(M)$ via the correspondence $b\to \omega_b$, provided that
$M$ is oriented.

\smallskip

{\bf Restricted transits.}  Pre-branched and branched moves are
naturally organized in subfamilies which give rise to restricted
transit equivalences.  In the branching setting, the so called ideal 
{\it sliding} transits  had been widely studied since
\cite{BP2}. We shortly
revisit these early results, with some complement, and outline 
a (partially conjectural) conceptually different approach to Theorem \ref{bI}
and eventually also to Theorem \ref{nakedI}. 
The basic idea is to point out some structures of differential topological
nature on $M$ which are carried by every branched ideal triangulation
$(T,b)$ of $\hat M$, are preserved by the sliding transits and can be 
modified by the full branched transits. Complements mainly
concern the so called carried {\it horizontal foliation}.
In the pre-branching setting, and assuming $M$ oriented,
the {\it non ambiguous} equivalence has been introduced  
in \cite{NA} and leads to interesting examples
of so called {\it non ambiguous structures}, also related to the theory of {\it taut} 
triangulations \cite{La}, with the associated {\it
reduced} quantum hyperbolic invariants (see also \cite{QT}).  The
non ambiguos relation lifts to a specialization of the sliding one 
on branched triangulations; we
point out a few specific insights into the branching theory, again in terms
of carried structures preserved by the non ambiguous
and which can be modified by the whole slidings transits. 
In \cite{2D} we have developed a 2D
counterpart of the present paper. This theme also emerged in
\cite{NA}, Section $5$ within a so called ``holographic'' approach to
3D non ambiguous structures.

\section{Generalities on decorated triangulations and their transits}
\label{generality}
\begin{figure}[ht]
\begin{center}
 \includegraphics[width=7cm]{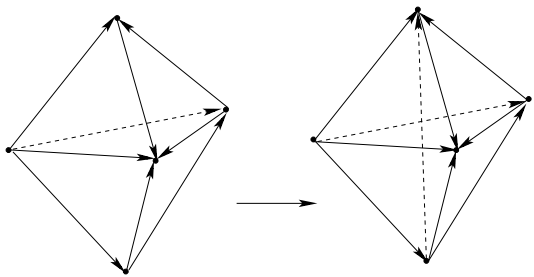}
\caption{\label{2-3} The $2\leftrightarrow 3$ move.}
\end{center}
\end{figure}

First let us recall the {\it ideal} naked moves. The $2\leftrightarrow 3$
move is illustrated in Figure \ref{2-3}; let us forget the edge
orientations to see the naked move; the arrows gives us also an
instance of branched move (see below). The positive move applies at
every couple of (abstract) tetrahedra with a common $2$-face which are
replaced by three tetrahedra around a new edge.  Both triangulations
share $6$ triangular faces on the boundary of the triangulated region.

\begin{figure}[ht]
\begin{center}
 \includegraphics[width=7cm]{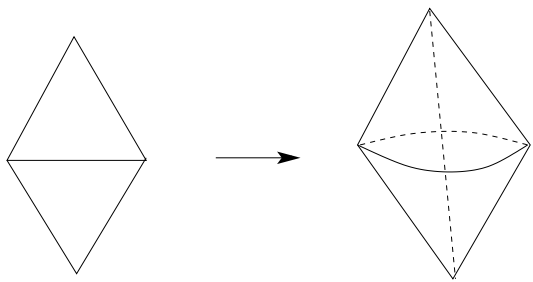}
\caption{\label{lune} The ideal quadrilater $0\leftrightarrow 2$ move.}
\end{center}
\end{figure}

Figure \ref{lune} illustrates the naked {\it quadrilateral}
$0\leftrightarrow 2$ ideal move. The positive one applies at every
quadrilateral made by two triangles with a common edge (a diagonal).
It is replaced by a ``pillow'' triangulated by two tetrahedra glued
along a quadrilateral which is a copy of the initial one on which a
diagonal exhange has been performed. The boundary of the pillow is
made by two copies of the initial quadrilateral glued along the common
boundary.

\begin{figure}[ht]
\begin{center}
 \includegraphics[width=8cm]{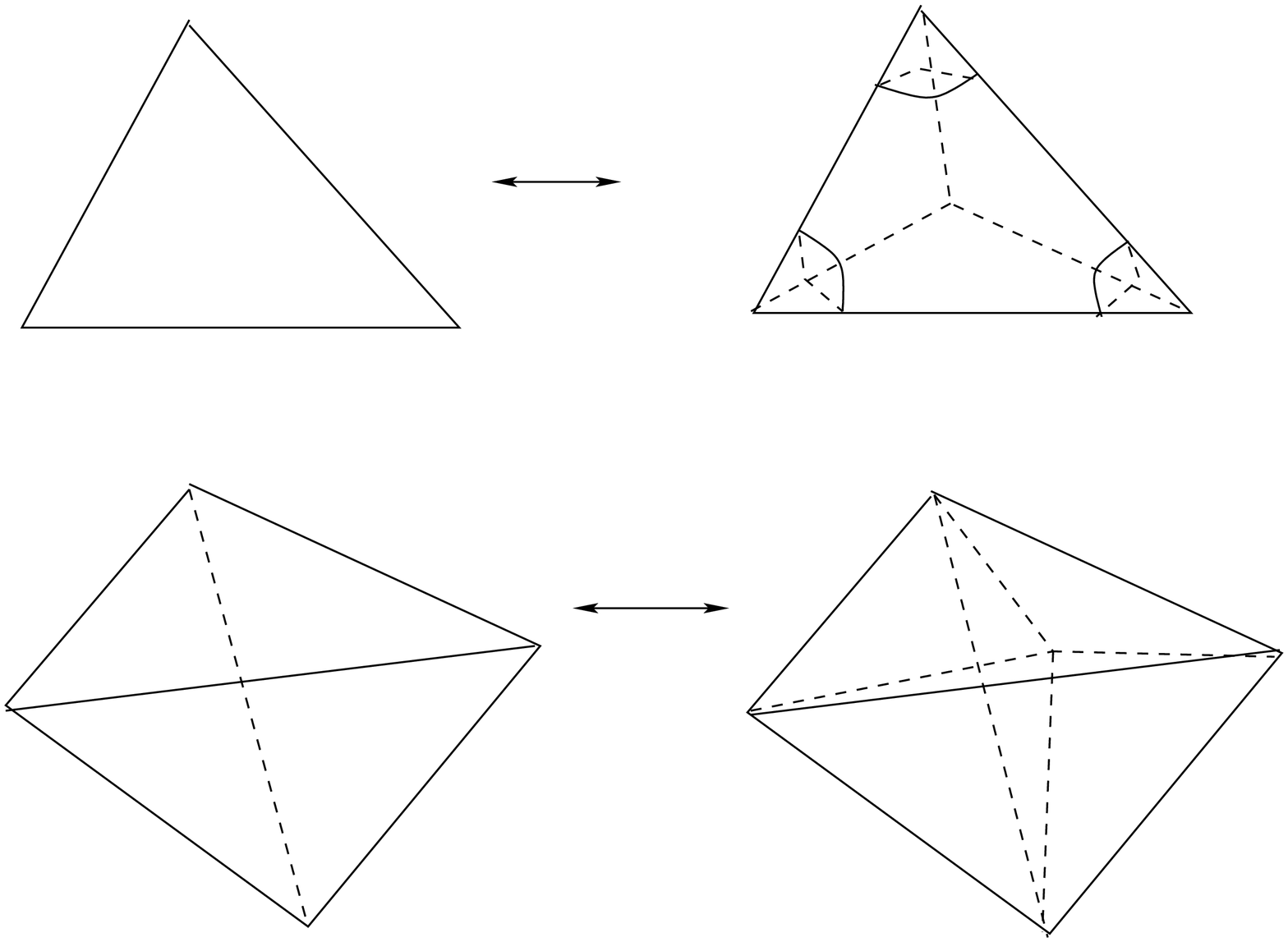}
 \caption{\label{completed} The triangular $0\leftrightarrow 2$ and
   the $1\leftrightarrow 4$ moves generating the completed transit
   relation.}
\end{center}
\end{figure}

\begin{figure}[ht]
\begin{center}
 \includegraphics[width=8cm]{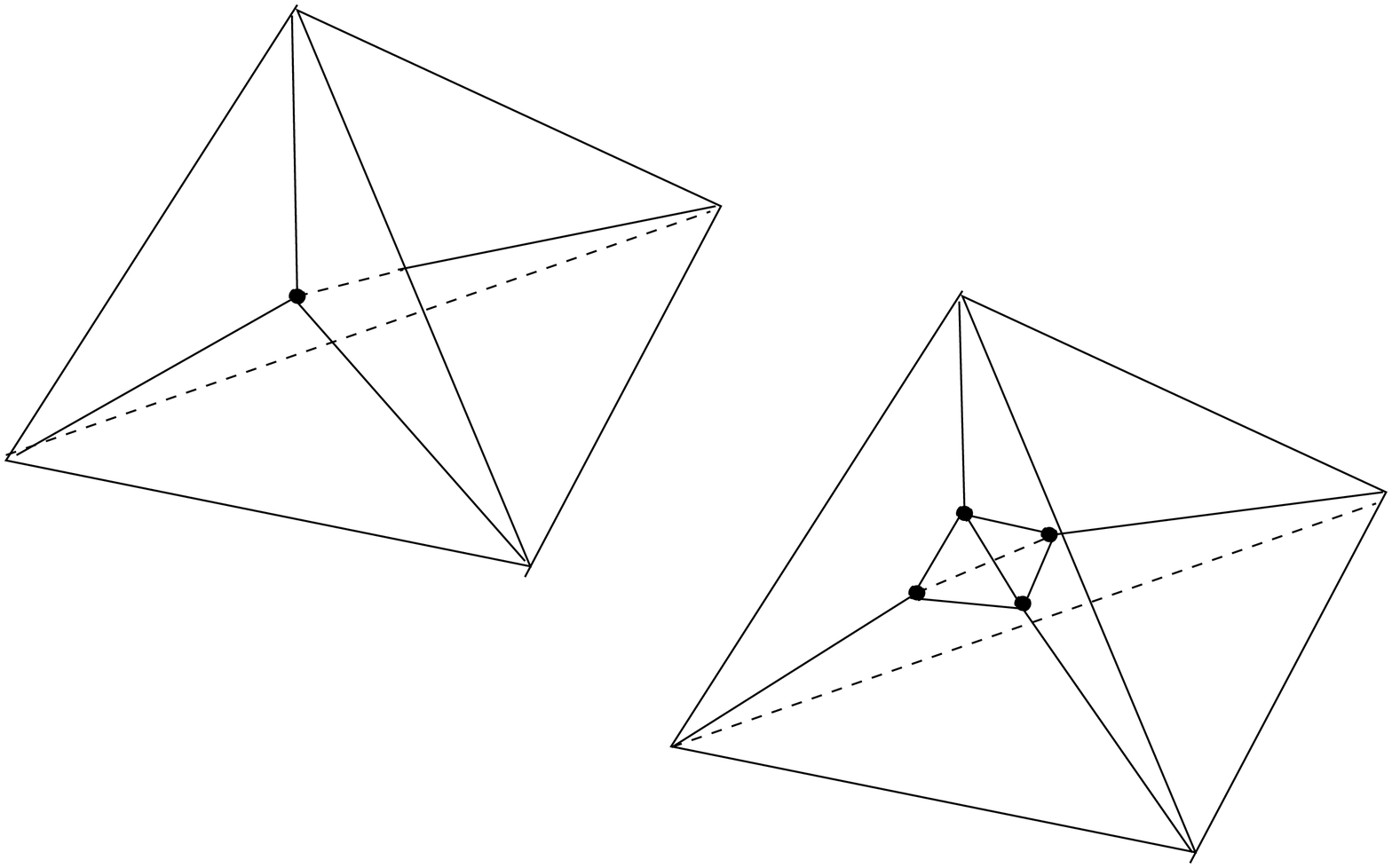}
 \caption{\label{butterfly} A butterfly and its modification by a stellar $1\to 4$ move.}
\end{center}
\end{figure}

In Figure \ref{completed} we show the two moves that generate
(together with the above ideal ones) the {\it completed} transit
relation. In fact we can eliminate one of them without modifying the
relation (see also \cite{RST}). We will freely use both. The positive
{\it triangular} $0\rightarrow 2$ move applies at every triangle. It is
repalced by a pillow triangulated by two tetrahedra glued along a copy
of the initial triangle on which a 2D stellar $1\to 3$ move has been
performed.  The boundary of the pillow is made by two copies of the
initial triangle glued along the common boundary. The positive $1\to 4$ move
is shown on the bottom of Figure \ref{completed}; in Figure \ref{butterfly}
we see how it acts on a dual butterfly.  

\smallskip

{\bf Decorated transit.} Let $T\to T'$ be a {\it positive} ideal naked
move. We describe the decorated enhancements $(T, d)\to (T', d')$
where $d, d'$ are either branchings or pre-branchings.  Consider first
a $2\to 3$ move. Then $T$ and $T'$ share (on the portion involved in
the move) 9 edges and 6 $2$-faces. If $d=b, d'=b'$ are branchings we
require that $b$ and $b'$ agree at these persistent 9 edges. If
$d=\omega, d'=\omega'$ are pre-branchings then we require that the
transverse orientations agree at those $6$ persistent $2$-faces. If we
consider a quadrilateral $0\to 2$ move, in the branching case we
require that the pillow boundary is made by two {\it branched} copies
of the initial branched triangulated quadrilateral. Similarly in the
case of pre-branchings we require that the transverse orientations at
the pillow boundary are two copies of the transversal orientations at
the triangles of the initial quadrilateral. A {\it negative}
$d$-transit is by definition the inverse of a positive one.  

If $M$ is {\it oriented}, every
branching $b$ incorporates a pre-branching $\omega_b$.  If $b\to b'$
is a $b$-transit, then it induces a $pb$-transit $\omega_b \to
\omega_{b'}$. It is easy to see that for every positive ideal transit
$T\to T'$ and for every decoration $(T,d)$ then we can enhace it to a
$d$-transit $(T,d)\to (T,d')$ (see below for more details).
On the other hand {\it there are decorations
$(T',d')$ such that the negative ideal move $T \leftarrow T'$ cannot
be enhaced to any $d$-transit}. Moreover, this happens for $(T',b')$ if
and only if it happens for $(T', \omega_{b'})$. Hence there is a
natural forgetting map
$$ \phi: \Bb^{id}(M) \to \Pp\Bb^{id}(M), \ \phi([T,b]) =
[(T,\omega_b)] \ . $$  

For the further two moves generating the completed relation we do
similarly.
\medskip

Here is a first useful reduction in order to study the decorated
transit equivalences. In fact next Lemma reduces the study
to decorations carried by a same triangulation $T$.

\begin{lem}\label{sameT} 
  For every $M$, for every triangulations $(T_1,d_1)$ and
  $(T_2, d_2)$ of $\hat M$,
  there are $(T,d)$ and $(T,d')$
  such that $(T_1,d_1)$ is equivalent to $(T,d)$ and
  $(T_2,d_2)$ is equivalent to $(T,d')$.
\end{lem}

\Dim By \cite{Ma}, there are finite sequences
of naked positive ideal moves $T_j \Rightarrow T$, $j=1,2$ (including both
$2\to 3$ and quadrilateral $0\to 2$ moves).  By
positivity, there is no obstruction to enhance them to sequences of
decorated transits.

\cvd

\medskip

\subsection{A combinatorial classification of decorated transits}\label{sub-fam}
\begin{defi}\label{forced}{\rm Given a {\it positive}
ideal move $T\to T'$ and a decoration $(T,d)$ we say that a
$d$-transit $(T,d)\to (T',d')$ is {\it forced} if it is the unique
enhancement of the naked move with the given initial decoration.}
\end{defi}
\smallskip

It is easy to see that if a decorated transit is not forced, then there are
exactly two such enhancements.

Let us fix some notations. A positive
$2\to 3$ move acts on the union of two ``abstract'' tetrahedra
$\tau_1\cup \tau_2$ which share one triangle $t$. We denote by $v_j\in
\tau_j$ the vertex opposite to $t$.  A naked positive quadrilateral
$0\to 2$ move acts on the union of two triangles $t_1\cup t_2$ which
share one edge $e$. We denote by $v_j\in t_j$ the vertex opposite to
$e$. 

\begin{defi}\label{sliding} {\rm A positive ideal $b$-transit
$(T,b)\to (T',b')$ is called a {\it bump transit} if the vertices $v_1$
    and $v_2$ (in $(\tau_j,b)$ or in $(t_j,b)$ as above) are both
    either a pit or a source.  A positive ideal $b$-transit $(T,b)\to
    (T',b')$ is called a {\it sliding transit} if it is not a bump transit.  
    A negative $b$-transit is either
    bump or sliding if it is the inverse of a positive one.}
\end{defi} 

We denote by $\Ss^{id}(M)$ the quotient set for the ideal $s$-transit
relation generated by sliding $b$-transits. Hence we have a
natural surjective projection
$$ \Ss^{id}(M)\to \Bb^{id}(M) \ . $$

Now we recall the notion of {\it non ambiguous transit}; in doing it  
we stipulate that $M$ is  {\it oriented}.

\begin{defi}\label{NAdef} {\rm 
(1) A positive ideal $pb$-transit 
$(T,\omega)\to (T',\omega')$ is
{\it non ambiguous} if it is forced.

(2) A positive ideal $b$-transit $(T,b)\to (T',b')$ is {\it non
  ambiguous} if $(T,\omega_b) \to (T',\omega_{b'})$ is non ambiguous.

(3) A negative decorated transit is non ambiguous if it is the
inverse of a non ambiguous positive one.}
  \end{defi}
  
  We stress that (2) is not the immediate definition of non ambiguous
  $b$-transit one would wonder. It is actually stronger. In fact we
  have (see below for more details):
  
  \begin{itemize}
  \item If $(T,b)\to (T,b')$ is non ambiguous, then it is forced.
  
  \item There are instances of forced $(T,b)\to (T',b')$ 
  which are ambiguous.
  
 \end{itemize}  
  
\smallskip

The non ambiguous transits define restricted $na$-transit
equivalences.  We denote by $\Nn\Aa^{id}(M)$ the corresponding
quotient set in the case of pre-branchings; $\Bb\Nn\Aa^{id}(M)$ the
one for branchings. We have a natural surjective projection
$$ \Bb\Nn\Aa^{id}(M)\to \Ss^{id}(M) \ . $$
Non ambiguous ideal $pb$-transits are widely characterized and
discussed in \cite{NA}. 
Now we make a more accurate analysis of the distributions of ideal
$b$-transits in such sub-families. We have defined above the signs
$*_{(\Delta,b)}$ when $M$ is oriented. If $M$ is not necessarily
oriented (orientable) we can fix anyway a local auxiliary orientation
on the (abstract) portion involved by a given move, and the signs
$*_{(\Delta,b)}$ are defined as well.  The initial configurations of a
positive ideal $2\to 3$ $b$-moves are in bijection with the total
orderings of the 5 vertices of $\tau_1\cup \tau_2$, hence they are
$120$. Every such an ordering induces a total ordering of the vertices
of both $\tau_1$ and $\tau_2$ which are respectively encoded by means
of the labels $\{0,1,2,3\}$ and determines the branchings
$(\tau_j,b)$. We can organize such configurations by $40$ {\it types}
which are classified by the couples
$$ ((*_1,a_1),(*_2,a_2))\in (\{\pm 1\} \times \{0,1,2,3\})^2$$ where
$*_j$ is the $b$-sign of $(\tau_j,b)$, $a_j$ is the order number of
$v_j$ with respect to the vertex ordering on $\tau_j$. Each type
carries $3$ configurations which form the orbit of the group of cyclic
permutations of the vertices of $t$.  The (non)ambiguous/sliding/bump
nature of a configuration only depends on its type; so the types are
distributed according to this classification.  Below we list $20$
types arranged along a few rows; each type has an implicit {\it
  symmetric} one obtained by exchange the two entries of the couple.
This corresponds to exchanging the role of $v_1$ and $v_2$. It is
understood that symmetric types belongs to the same row. So we have:

\begin{itemize}

\item {\bf Non ambiguous $2\to 3$ $b$-transits.}
$$((-1,1), (-1,0)),\  ((+1,1), (+1,0)),\  ((+1,2),(+1,3)),\   ((-1,2),(-1,3))$$
$$ ((+1,2),(-1,0)),\ ((-1,3),(+1,1)), \ ((-1,2),(+1,0)),\ ((+1,3),(-1,1))$$
$$((+1,2),(+1,1)), \ ((-1,2),(-1,1)) \ . $$
The last two are characterized by the fact that all $5$ tetrahedra
occurring in the transit share the same sign $\pm 1$ (the example shown
in Figure \ref{2-3} is of this kind); 
sometimes we call them {\it Schaeffer's}  transits and they play a distinguished
role in the study of quantum hyperbolic state sums (see \cite{NA}, \cite{AGT}).

\medskip

\item {\bf Ambiguous sliding $2\to 3$ $b$-transits.}
$$((-1,1),(+1,1)), \ ((+1,1),(-1,1)), \ ((+1,2),(-1,2)), \ ((-1,2),(+1,2))$$
$$((-1,3),(-1,0)), \ ((+1,3),(+1,0)) \ . $$
The last two are characterized by the fact that although they are ambiguous
(at the level of the induced pre-branchings), they are forced;
so they are
called {\it forced ambiguous $b$-transits}.

\medskip

\item{\bf $2\to 3$ bump transits.}
$$((+1,0),(-1,0)), \ ((-1,0),(+1,0)), \ 
((+1,3),(-1,3)) \ ((-1,3),(+1,3)) \ . $$

\end{itemize}

\medskip

Note that row by row, there are pairs of cases
just related by the inversion of the $b$-signs.  This reflects
geometric symmetries and instances in a pair basically have the same
qualitative features.  A very similar combinatorial classification
(including the presence of forced ambiguous instances) holds for the
quadrilateral $0\to 2$ move.

\begin{figure}[ht]
\begin{center}
 \includegraphics[width=12cm]{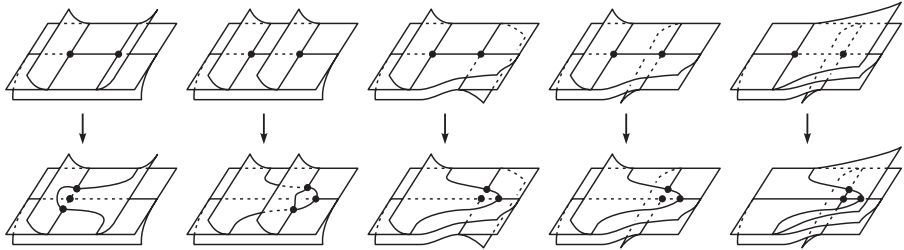}
 \caption{\label{MPSLIDE} The sliding branched spine $2\to 3$ move.}
\end{center}
\end{figure}

\begin{figure}[ht]
\begin{center}
 \includegraphics[width=8cm]{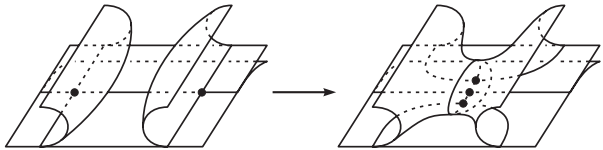}
 \caption{\label{MPbump} A branched spine $2\to 3$ bump move.}
\end{center}
\end{figure}

\begin{figure}[ht]
\begin{center}
 \includegraphics[width=10cm]{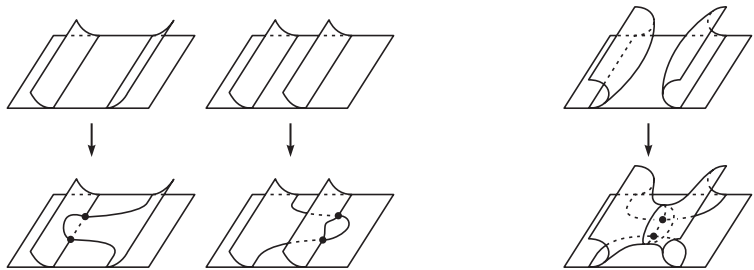}
 \caption{\label{lunebr} Branched lune moves.}
\end{center}
\end{figure}

\begin{remark}\label{why-slide}{\rm
    In Figures \ref{MPSLIDE}, \ref{MPbump}, \ref{lunebr} we see instances of the above classification
``sliding vs bump'' moves in terms of dual spines. A comment is in order.
    Every abstract tetrahedron and the dual butterfly have its own full group of symmetries.
    So a `correct' picture of a naked moves should reflect them.
For example in Figure 17 (B) in \cite{RST} we see such a symmetric picture of
the $2\to 3$ move in terms of spines. On the other hand, in Figure 18 of \cite{RST} we see another current
picture of it; here the symmetry  has been {\it broken} at both vertices of the initial spine and the moves 
becomes visually transparent: a region slides  across a vertex to create two new ones.
However this is somehow a visual artefact because such a symmetry break has no intrinsic meaning
on a naked triangulation. Similarly in Figure \ref{butterfly} we show a symmetric picture of the $1\to 4$ move in terms
of spines, while in Figure 17 (A) of \cite{RST} we see a current non symmetric one.
A branching can be also considered as a way to break the symmetry at every tretrahedron of a triangulation
(see Figure \ref{b+}) under a certain global coherence. We realize that sliding transits (in terms
of dual branched spines) are realized by actual region sliding, while along a bump transit a region bumps into another one.     
We will return on this point in Section \ref{slide}.
}
\end{remark}

 \medskip

\section{The branched connectivity results}\label{b-connectivity}
Let $(T,b)$ be an ideal branched triangulation of $\hat M$.
\begin{defi}\label{amb-edge}{\rm

    (1) An edge of $T$ is said {\it ambiguous} if 
    we can invert its orientation and keep a branched triangulation $(T,b')$
    (clearly the edge is ambiguous also in $(T,b')$).    
  
   \smallskip

    (2) An ambiguous edge $e$ of $(T,b)$ is {\it good} if:
    \begin{itemize}
    \item (i) It does not support self-guing, that is it lifts to at most one edge
    in every abstract tetrahedron of $T$.
    
    \item (ii) Consider the  branched star of $e$ and its link $L(e)$ formed by the union
    of the opposite edge to $e$ in each tetrahedron containing $e$. Then the branching
    does {\it not} make $L(e)$ an oriented circle.
    \end{itemize}  
  }
\end{defi}

We have

\begin{lem} Let $(T,b)$ and $(T',b')$ differ by the inversion of one
  good ambiguous edge. Then they are ideally $b$-transit equivalent
  to each other. 
\end{lem}
\Dim Let $e$ be the good ambiguous edge which is inverted passing from $b$ to $b'$.
Then $e$ must be ambiguous in every branched tetrahedron of its star.
Notice that in a branched tetrahedron with ordered vertices $v_0,v_1,v_2,v_3$,
the ambiguous edges are $v_0v_1$, $v_1v_2$, $v_2v_3$.
If the star consists of two tetrahedra, we recognize
the final configuration of a positive ambiguous quadrilateral $0\to 2$
$b$-transit. So we get the inversion of $e$ by pre-composing this move
with its negative inverse. If the star has more than two tetrahedra, 
by performing a positive $2\to 3$
$b$-move at two tetrahedra in the star at which the two edges in $L(e)$ have conflicting
orientations, $e$ persists being good ambiguous and that number decreases by 1.

\cvd

Then we can add the elementary move of {\it inverting a good ambiguous edge}
without modifying the ideal $b$-transit equivalence.

\subsection{Proof of Theorem \ref{bC}}\label{proof-bC}  
Thanks to Lemma \ref{sameT}, it is enough to proof that for every ideal triangulation
$T$, for every branchings $b$ and $b'$, $[T,b]=[T,b']\in \Bb(M)$. We perform
the first
two steps of the construction of the (naked) barycentric subdivision of
$T$ by means of positive $2\to 3$ and $1\to 4$ moves (see \cite{RST}, Section 2.5).
First we apply the $1\to 4$ move at every tetrahedron of $T$. We get a triangulation
$T_1$ in which the $2$-faces of $T$ persist. Every such a face $F$ is a common face
of two tetrahedra of $T_1$; perform a further $1\to 4$ move at one of them. Do it for every
$F$ getting a triangulation $T_2$. The faces $F$ persit also in $T_2$. Every $F$ is a common 
face of two tetrahedra of $T_2$. Perfom the $2\to 3$ move that eliminate $F$.
Do it for every $F$ as above. We get our final triangulation $\tilde T$.
Now let us enhance the composite transit  $T\Rightarrow \tilde T$ to branched transits
$(T,b)\Rightarrow (\tilde T, \tilde b)$ and   $(T,b')\Rightarrow (\tilde T, \tilde b')$
as follows: all branched $1\to 4$ moves are such that the new vertex is a pit.
The final $2\to 3$ $b$-moves are bump moves and we stipulate that the new
edges created by these moves are oriented in both $(\tilde T, \tilde b)$
and $(\tilde T, \tilde b')$ in such a way that the vertices created by some $1\to 4$ move keep
the property of being a pit. The naked edges of $T$ persist in $\tilde T$, and $\tilde b$
can differ from $\tilde b'$ only at some of these persistent edges. Finally
we readily see that they are all good ambiguous in both $(\tilde T, \tilde b)$ 
and $(\tilde T, \tilde b')$,
hence we can pass from $\tilde b$ to $\tilde b'$ by inverting some of them.

\cvd
 
\subsection{Proof of Theorem \ref {bI}}\label{proof-bI}
Again by Lemma \ref{sameT}, it is enough to prove that
for every ideal triangulation $T$ of $\hat M$, for every branchings $b$ and $b'$,
there is a composite ideal $b$-transit $(T,b)\Rightarrow (T,b')$.
We start with the composite non ideal $b$-transit
$$ (T,b) \Rightarrow (T_1, b_1) \Rightarrow (\tilde T,\tilde b)\Rightarrow
(\tilde T,\tilde b') \Leftarrow (T_1,b'_1)\Leftarrow (T,b')$$
constructed in the above proof.
Now we want to apply branched  versions of Matveev's {\it arch} and associated 
constructions to every  $1\to 4$ moves in order to eventually
replace them in a systematic way with a 
composition of branched ideal moves. Naked arch and associated constructions
are clearly discussed in \cite{RST}, Section 5. By forgetting for a moment
the branchings (so that the inversion of ambiguous edges becomes immaterial), 
our starting non ideal composite transits are much simpler  
than the general ones occurring in the proof of Theorem 1.4 in \cite{RST}.
The main reason is that in our situation there are no interferences
between the ideal and non ideal moves occurring in the given composite
transit; then we can apply straightforwardly  the naked construction described in 
Section 6 of \cite {RST}
without any need of the  subtle  ``moving the arch"  procedure developed therein. 
Moreover we have the advantage of allowing the quadrilateral $0\leftrightarrow 2$ move. 
So we can try to enhace such rather simple construction in the branched setting. 
\smallskip

\begin{figure}[ht]
\begin{center}
 \includegraphics[width=6cm]{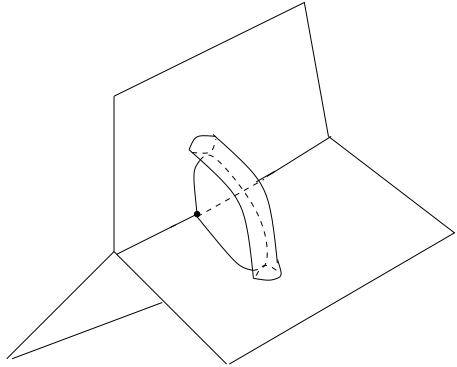}
 \caption{\label{arch} Building an arch at an edge of a spine.}
\end{center}
\end{figure}

First let us decribe the arch construction associated to every naked
$1\to 4$ move $\Tt_1 \to \Tt_2$ producing a new material vertex
$v$. One implements the following procedure.  After the move there are
six triangles $t$ with common vertex $v$, so that the edge of $t$
opposite to $v$ belongs to $\Tt_1$.  Select one $t$. Then select one,
say $e$, of the two edges of $t$ with vertex $v$; so in total there
are 12 possible choices of the couple $(t,e)$, each one called a {\it
  arch marking} (at the given naked $1\to 4$ move).  Every marking
encodes the location of an arch associated to the move.  This is a
quite localized construction that, roughly speaking, unglues the
triangulation $\Tt_2$ at $t$ an inserts one tetrahedron (with
face/edge identifications) in such a way that the two vertices of $e$
are eventually identified. It is easier to visualize an arch on the
dual spine.  The triangle $t$ corresponds to a transverse spine
edge. The arch opens a tunnel connecting the two boundary components
of $M$ corresponding to the vertices of $e$.  In Figure \ref{arch} we
show a realistic picture of it.  Assume now that the move
$(\Tt_1,\bG_1) \to (\Tt_2,\bG_2)$ is branched in such a way that $v$
is a pit.  We will deal only with this instance of branched non ideal
move so we will understand it.  The selected triangle $t$ is now
branched by $\bG_2$ and the edge $e'$ of $t$ opposite to $v$ keeps the
initial $\bG_1$-orientation. Then there are two possible choices for
the edge $e$ which are distinguished by the branching: either the
$\bG_2$-orientation of $e$ conflicts with the $\bG_1$-orientation of
$e'$ at the common vertex $e\cap e'$ or these orientations match; in
this second case we say that $(e,\bG_2)$ carries the {\it prevailing
  orientation}. We have:
\smallskip

{\bf Claim 1.} The naked arch encoded by the marking $(t, e)$ can be
enhaced to a branched arch (producing a branched triangulation
$\Tt'_2,\bG'_2)$ which modifies only locally $(\Tt_2,\bG_2)$) if and
only if $(e,\bG_2)$ carries the prevailing orientation. Hence there
are in total 6 possible implementations of the branched arch
construction associated to our branched $1\to 4$ move.
\medskip

{\bf Claim 2.} The composition of a branched $1\to 4$ move
followed by the insertion of an associated  branched arch
can be undone (reobtaining the initial branched tertahedron $(\Delta,b)$ supporting
 the $1\to 4$ move) by means of ideal $b$-moves localized on $\Delta$. 
\medskip

First we complete the proof by assuming the Claims. Finally we will prove them.
\smallskip

A few remarks about our starting composite transit:
\begin{enumerate}
\item Forgetting the branchings, the left and right composite transits
are two copies of $ T \Rightarrow T_1 \Rightarrow \tilde T$.

\item On the left side, every edge which is present at some stage
of the composite transit $(T,b)\Rightarrow (\tilde T, \tilde b')$
persits for ever keeping its branching orientation. The same
fact holds on the right side composite transit $(\tilde T, \tilde b') \Leftarrow (T,b')$.
\end{enumerate}

So every branched non ideal move $(\Tt_1,\bG_1)\to (\Tt_2,\bG_2)$
occurring on the left side has a twin move $(\Tt_2,\bG'_2)\leftarrow
(\Tt_1,\bG'_1)$ on the right.  Now we have to specify the (branched)
arch marking at every $1\to 4$ moves.  For the move occurring in
$$(T_1, b_1) \Rightarrow (\tilde T,\tilde b)\Rightarrow (\tilde
T,\tilde b') \Leftarrow (T_1,b'_1)$$ we can choose two copies of the
same $(t,e)$ in every couple of twin moves, such that
$(e,\bG_2)=(e,\bG'_2)$, both carry the prevailing orientation, and $e$
is not a persistent edge coming from $T$.  For the couples of twin
moves occurring in
$$ (T,b) \Rightarrow (T_1, b_1), \ (T_1,b'_1)\Leftarrow (T,b')$$ we
take two copies of the same $t$ and two copies of the same $e$
(carrying the prevailing orientation) if the $\bG_1$ and $\bG'_1$
orientation agree on $e'$; otherwise the two branched arch markings
differ by the choice of the edge $e$.  Now one readily checks that no
triangle occurring in a marking is destroyed by any ideal $b$-move of
the composite transits. Then we can conclude by applying the above
Claims, formally in the same way that works in the naked proof
(\cite{RST}) Section 6) in the very favourable circumstance that
``moving the marking" is never necessary.  It remains to justify the
Claims.
\smallskip 

{\it Claim 1:} Again it is easier to deal in terms of dual spine.
Possibly by using an auxiliary local orientation of $M$ we can assume
that the three involved portions of spine regions are oriented
satisfying the branching condition. To build the arch we remove a
small open $2$-disk from two of them and we attach a copy of
$S^1\times [0,1]$ glued to a monogon $\tau$ in such a way that
$(S^1\times [0,1])\cap \tau= \{p\}\times [0,1]$.  We need that the
given region orientations propagate across $(S^1\times [0,1])\setminus
(\{p\}\times [0,1])$ in such a way that the branching condition is
kept; note that the monogon, that is its dual edge in
$(\Tt'_2,\bG'_2)$, is ambiguous. It is now easy to check that this
happens just when $e$ carries the prevailing orientation.
For completeness we describe also the
corresponding one branched tetrahedron of $(\Tt'_2,\bG'_2)$ (with
face/edge identification) inserted at the triangle $t$ to build the
arch.  By adopting the notations stated above, let $v_0,v_1,v_2$ be
the ordered vertices of $(t,\bG_2)$, so that $v=v_2$, $e= v_1v_2$.
Realize $(t,\bG_2)$ as the $2$-face of an abstract branched
tetrahedron $(\Delta, \Bb)$ in such a way that its opposite vertex $w$
is smaller that both $v_2$ and $v_1$. Note that this requirement makes
ambiguous the edge of $\Delta$ with vertices $w$ and $v_0$. Now let us
identify two $2$-faces of $(\Delta,\Bb)$ respecting the branching;
precisely we identify $(t,\bG_2)$, that is the $2$-face $v_0,v_1,v_2$,
with the $2$-face $w,v_1,v_2$, so that $w$ and $v_0$ are identified at
a same point $p$ . Finally let us identify the edges $pv_1$ and
$pv_2$, so that $v_1$ is identified with $v_2$.

\cvd

\begin{figure}[ht]
\begin{center}
 \includegraphics[width=8cm]{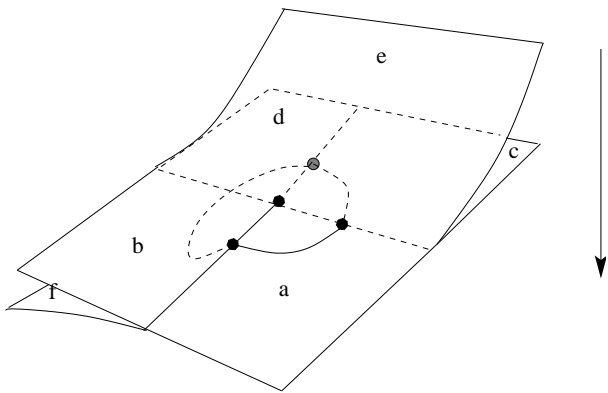}
 \caption{\label{b-arch1} A branched $1\to 4$ move at a branched butterfly.}
\end{center}
\end{figure}

{\it Claim 2:} It is enough to follow {\it backward} the sequence of
spine moves given in Figure 9 of \cite{Am} (or equivalently the same
sequence splitted in two parts accordingly to Figure 29 and 28 of
\cite{RST}) by checking that it can be performed in the branched
setting. Let us start with a branched butterfly as in Figure \ref{b-arch1}.

\begin{figure}[ht]
\begin{center}
 \includegraphics[width=8cm]{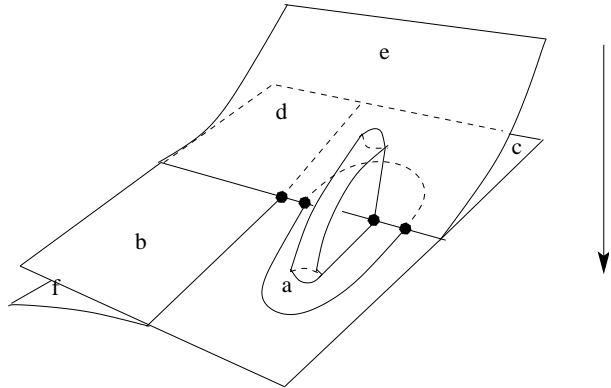}
 \caption{\label{b-arch2} A $b$-arch at a branched bubble move.}
\end{center}
\end{figure}

We assume that the transverse orientation at everyone of the $6$ germs
of region of the butterfly is given by the vertical vector pointing
towards the bottom. In the picture these germs are labelled by $a$,
$b$, $c$, $d$, $e$, $f$. Strictly speaking there should be a few such
configurations to consider, also taking into account the (local) sign
of the dual branched tetrahedron. We limit to discuss one as the other
cases do not present substantial differences.  The branched $1\to 4$
move is obtained dually by attanching a $2$-disk $D$ along an embedded
smooth circle into the smooth sheet of the butterfly made by the union
of the regions labelled by $a$, $c$, $f$.  This forms a bubble whose
interior contains the new material vertex $v$. Moreover $D$ is
transversely oriented by a vertical vector pointing towards the top
(accordingly to the fact that $v$ is a pit in the branched
triangulation). A branched $1\to 4$ move, expressed in spine terms, is
the composition of a branched bubble move followed by a positive
branched $2\to 3$ move.  If we had built a (suitable) $b$-arch at the
$1\to 4$ move, this follows the inverse $3\to 2$ move, so that we
reach a configuration of $b$-arch at a branched bubble move as it is
illustrated in Figure \ref{b-arch2}. The $2$-disk $D$ is now glued
along a smooth circle embedded in the sheet made by the union of the
regions labelled by $a$ and $c$.

\begin{figure}[ht]
\begin{center}
 \includegraphics[width=8cm]{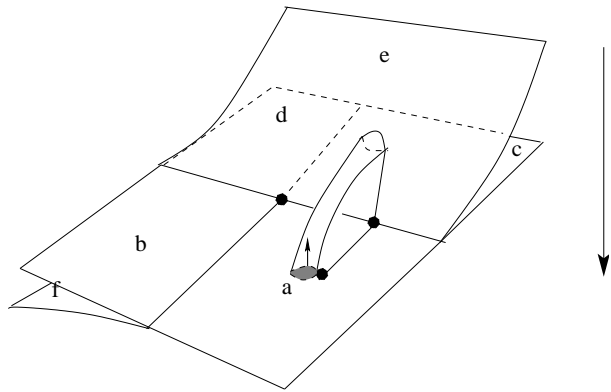}
 \caption{\label{b-arch3} A $b$-arch with membrane.}
\end{center}
\end{figure}

Now the disk $D$ slides by means of a further branched $3\to 2$ move
realizing (by using the terminology of \cite{RST}) a {\it $b$-arch
  with membrane} as shown in Figure \ref{b-arch3}.The picture also
shows the tranverse orientation of this membrane.

\begin{figure}[ht]
\begin{center}
 \includegraphics[width=8cm]{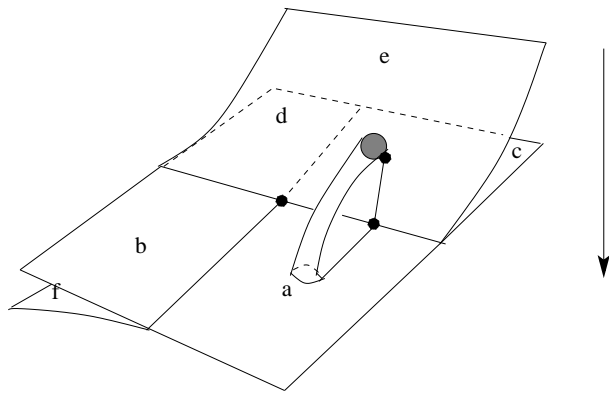}
 \caption{\label{b-arch4} A  $b$-arch with membrane again.}
\end{center}
\end{figure}

We can move by isotopy the membrane to reach the configuration of
Figure \ref{b-arch4}. Now the transverse orientation of the membrane
coincides with the one of the region labelled by $e$.

\begin{figure}[ht]
\begin{center}
 \includegraphics[width=8cm]{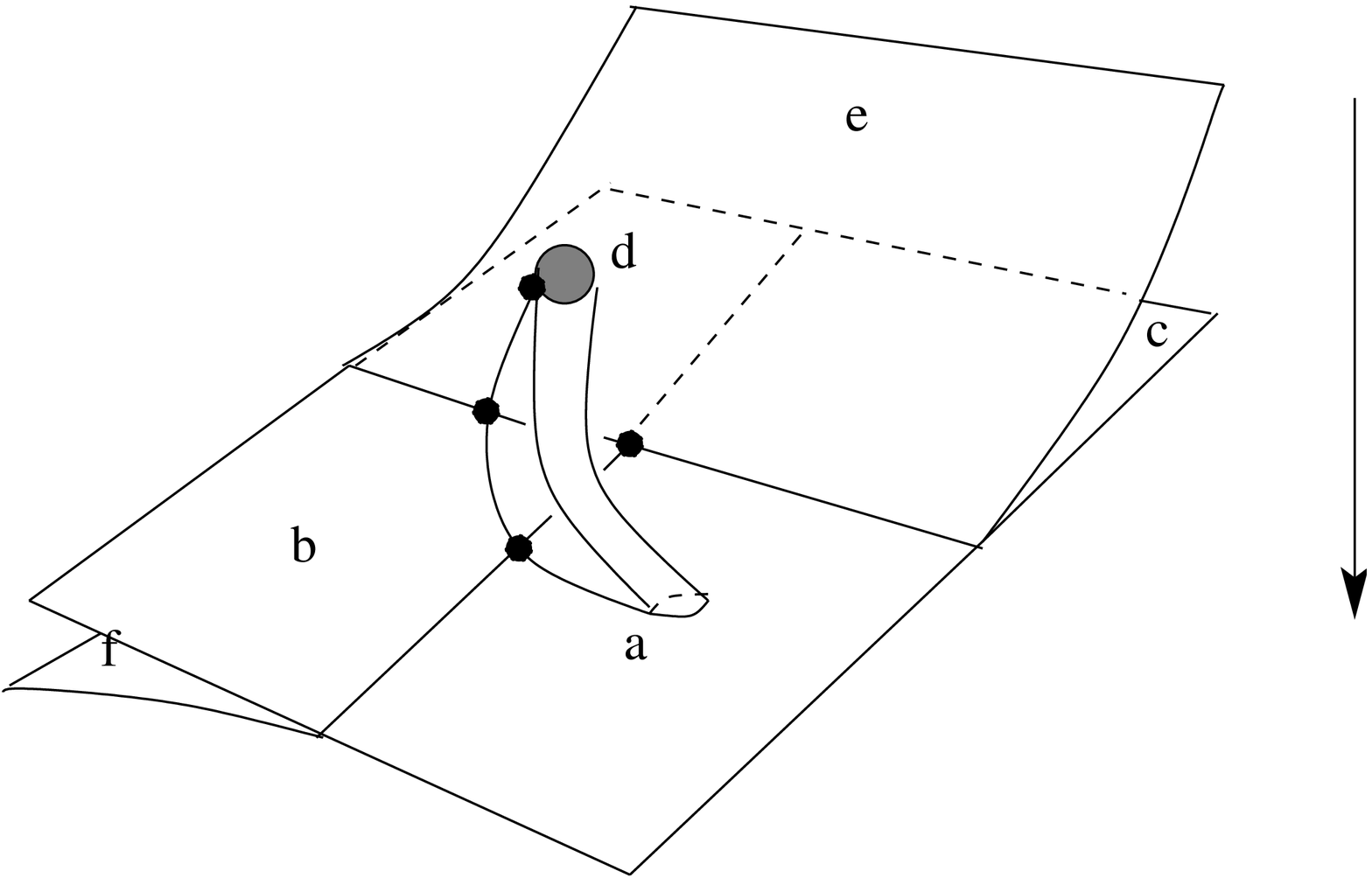}
 \caption{\label{b-arch5} Sliding the $b$-arch with membrane.}
\end{center}
\end{figure}

Now by means of a branched $2\to 3$ move we reach the configuration
illustrated in Figure \ref{b-arch5}.

\ref{b-arch5}

\ref{b-arch5}
\begin{figure}[ht]
\begin{center}
 \includegraphics[width=8cm]{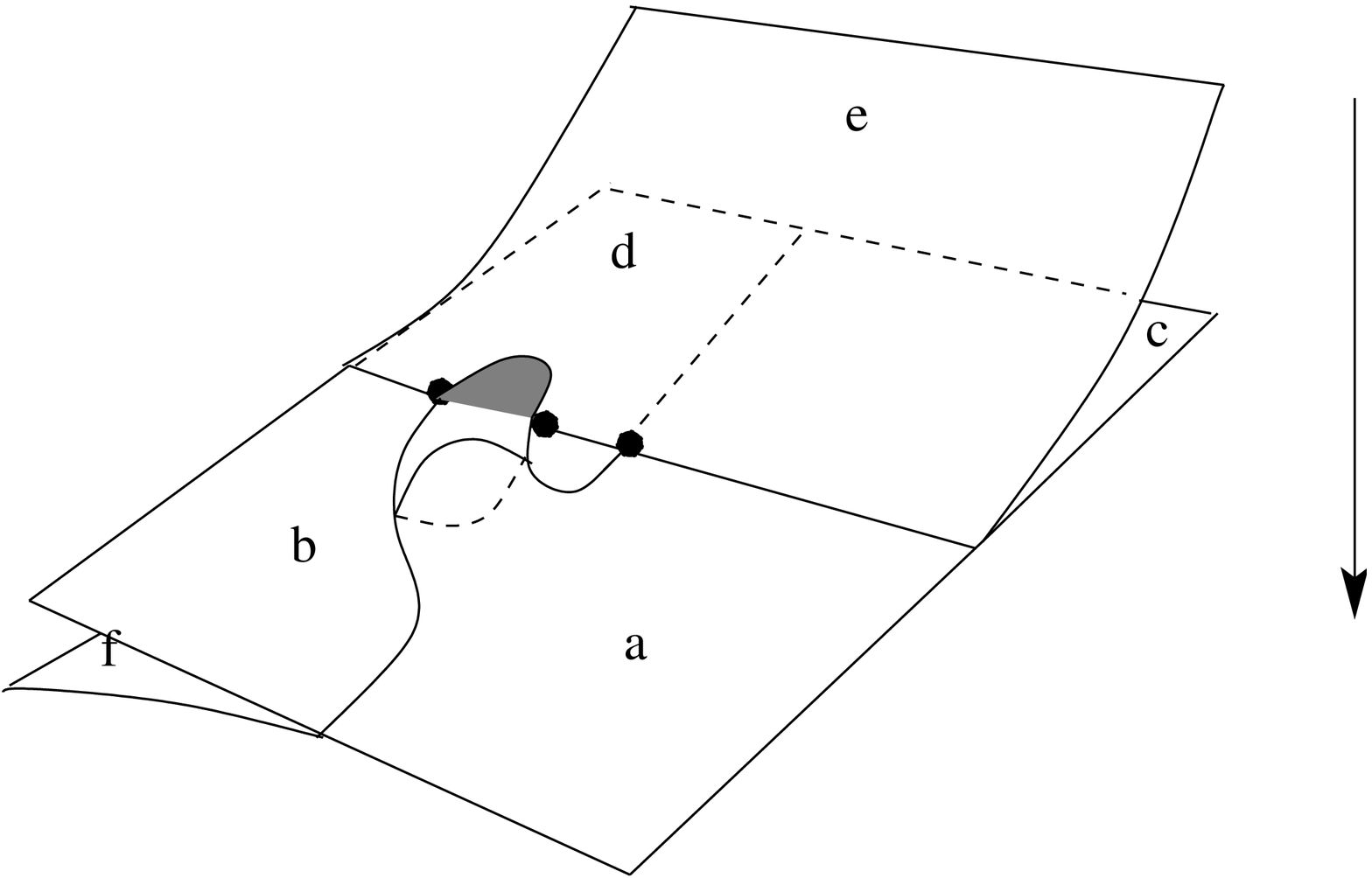}
 \caption{\label{b-arch6} A positive sliding lune move.}
\end{center}
\end{figure}

By performing a further branched $2\to 3$ move we reach the
configuration illustrated in Figure \ref{b-arch6}. Finally this can be
undone (getting the initial branched butterfly) by means of a negative
sliding lune move. Both Claims, hence Theorem \ref{bI} are eventually
achieved.

\cvd

\section{On $\Pp\Bb^{id}(*)$}\label{on-pb-b}
 As for pre-branching we mainly refer to \cite{NA}, accordingly with 
 this reference  in this section $M$ is assumed {\it oriented}.

\subsection{ Homological invariants}\label{homol-Inv}  
Let $T$ and $\Sigma$ be as usual.  Every pre-branching
$(T,\omega)$ can be interpreted as a {\it fundamental} cellular
$1$-$\Z$-cycle supported by the whole $1$-skeleton Sing$(\Sigma)$ of
$\Sigma$.  As a consequence of the ``circulation lemma'' of \cite{GR},
a branching $b$ can be interpreted as a {\it fundamental} cellular
$2$-$\Z$-chain supported by the whole $2$-skeleton of $\Sigma$ such
that the boundary $\partial b$ is a pre-branching (so that
$\omega_b=\partial b$).  Here ``fundamental'' means that all chain
coefficients are equal to $\pm 1$ (with respect to an arbitrary
auxiliary system of orientations of the cells of $\Sigma$). Given a
branching $b$, $\omega_b=\partial b$, we can take these orientations
in such a way that all coefficients are equal to $1$. 
We say that a
cellular chain is {\it almost fundamental} if the coefficients belong
to $\{0,\pm 1\}$. We have
 
 \begin{lem}\label{h-inv}  For every pre-branching $(T,\omega)$,
 
 (1) The class $[\omega]\in H_1(M;\Z)$ is invariant under
 ideal $pb$-transit equivalence.

(2) The reduction mod$(2)$ $[\omega]_2=0 \in H_1(M;\Z/2\Z)$.
 
(3) The class $[\omega]$ is even, i.e. $[\omega]=2\alpha$
 for some $\alpha \in H_1(M;\Z)$.
 
(4) If $\omega=\omega_b$ for some branching $b$, then $[\omega]=0$. 
 \end{lem}
 
 \Dim Point (1) follows immediately by looking at the local transits. 
 Point (2) holds because (forgetting the orientation)
 $\omega$ is the boundary of the unique fundamental $\Z/2\Z$-$2$-chain
 on $\Sigma$. Point (4) has been remarked above. As for (3), it is a
 general fact that $[\omega]$ is even if and only if its reduction
 mod$(2)$ $[\omega]_2= 0 \in H_1(M;\Z/2\Z)$.  Let us give anyway a
 constructive proof in the spirit of transit equivalence.
 Assume that $(T,\omega)$
 is such that $T$ carries some branching $b$. Then $[\omega]-[\partial
 b]=[\omega]=[\omega - \partial b]$ and this last is a cycle of the
form $z=2a$, where $a$ is an almost fundamental cycle. Then
$[\omega]=2[a]$. In general, in \cite{BP2} Theorem
 3.4.9 one finds an algorithm that for every naked triangulation $T$
 of $\bar M$, produces a chain of positive $2\to 3$ moves
 $T\Rightarrow T'$ such that $T'$ carries some branching $b$. Enhace
 this to $(T,\omega)\Rightarrow (T',\omega')$.  Then
 $[\omega]=[\omega']$ and the above argument applies to $\omega'$ and
 $b$.
 
 \cvd 
 
 \subsection{Non triviality of $\Pp\Bb^{id}(*)$}
 By Theorem \ref{bI}, the image of the forgetting map
 $$\phi: \Bb^{id}(M)\to \Pp\Bb^{id}(M)$$ 
 consists of one point. Without using this result, we are going to see anyway
that in general this image is a proper subset of $\Pp\Bb^{id}(M)$. 
We will see later that $\Pp\Bb^{id}(M)$ can be even infinite.  

\begin{lem}\label{Imphi} $[(T,\omega)]$ belongs to the image of
  $\phi: \Bb^{id}(M)\to \Pp\Bb^{id}(M)$ if and only if there exists a
  branching $(T,b)$ such that $[(T,\omega)]=[(T,\omega_b)] \in
  \Pp\Bb^{id}(M)$.
\end{lem}
\Dim ``If'' is trivial. On the other hand, let $(T,\omega) \Rightarrow
(T', \omega_{b'})$ a chain of $pb$-ideal moves such that $(T',b')$ is
branched. We can enhance the inverse naked chain $T\Leftarrow T'$ to a
$b$-chain $(T,b)\Leftarrow (T',b')$; in fact there is not any stop
because there is not at the pre-branching level. Finally
$[(T,\omega)]=[(T,\omega_b)]$.

\cvd

\begin{remark}{\rm By the above proof we cannot conclude that $\omega
    = \omega_b$ (for some suitable implementation of $(T,b)\Leftarrow
    (T',b')$) because of the possible presence of forced ambiguous
    $b$-transit which might change the pre-branching at some step.}
\end{remark}

We know that a necessary condition in order that $(T,\omega)$
represents a point in Im$(\phi)$ is that $[\omega]=0 \in
H_1(M;\Z)$. In some case this is also sufficient; forthcoming
  examples will show that the hypothesis of next Lemma \ref{H2=0} is
  sharp.

 \begin{lem}\label{H2=0} 
  Assume that $H_2(M;\Z/2\Z)=0$. 

  (1) Let $(T,\omega)$ be a pre-branched triangulation of $\hat M$
  such that $[\omega]=0\in H_1(M;\Z)$.  Then there exists a branching
  $(T,b)$ such that $\omega=\omega_b$.  

  (2) A triangulation $(T,\omega)$ represents a point in Im$(\phi)$ if
  and only if $[\omega]=0\in H_1(M;\Z)$.
\end{lem}
\Dim Clearly $(1)\Rightarrow (2)$. As for (1), let
 $$\omega = \partial \beta, \ \beta = 
 \sum_{C\in \Sigma^{(2)}} \beta(C)C $$ where $\beta$ is a cellular
 (not necessarily fundamental) $\Z$-$2$-chain on $\Sigma$.  Consider
 the union of $2$-cells such that $\beta(C)$ is even. This is a
 cellular $\Z/2\Z$-cycle of $\Sigma$.  As $\Sigma$ is a spine of $M$
 and $H_2(M;\Z/2\Z)=0$, then this is necessarily empty.  It follows
 that for every $2$-cell $C$, $\beta(C)$ is odd. Define $b= \sum_C
 \frac{\beta(C)}{|\beta(C)|}C$.  This is a branching and $\partial b =
 \omega$.

\cvd

\medskip

In some case a branching is determined by the associated
pre-branching.

\begin{lem}\label{H2Z=0} Assume that $H_2(M;\Z)=0$. 
  Let $(T,b)$ and $(T,b')$ be branched triangulations of $\hat M$ such
  that $\omega_b=\omega_{b'}$, then $b=b'$.
\end{lem}  
\Dim Let us orient the $1$- and $2$-cells of $\Sigma$ by $\omega_b$
and $b$ respectively, so that both $\omega_b$ and $b$ considered as
fundamental $\Z$-chains have all coefficients equal to $1$.  Set
$\beta=b-b'$, $\partial \beta = \omega_b - \omega_b'=\alpha$.  For
every $2$-cell $C$, $\beta(C)\in \{0,2\}$ as well as for every edge
$e\subset {\rm Sing}(\Sigma)$, $\alpha(e)\in \{0,2\}$.  Hence,
$\partial \frac{\beta}{2}=\frac{\alpha}{2}$.  As
$\omega_b=\omega_{b'}$, then $\frac{\beta}{2}$ is an almost
fundamental $\Z$-$2$-cycle which is empty because $H_2(M;\Z)=0$.

\cvd

\begin{remark}\label{H2Z-sharp} {\rm The hypothesis of Lemma
  \ref{H2Z=0} is sharp. In fact, $H_2(M;\Z)$ is free; if it is not
  zero, according to the above proof, there exists $(T,b)$ such that
  $b$ can be modified to some $b'= b-\partial 2z$, where $z$ has
  coefficients in $\{0, \pm1\}$, $b\neq b'$, and $\omega_b=\omega_{b'}$ (see also the
  proof of Proposition \ref{infinite-pb}).}
\end{remark}     
 
\medskip

{\bf Claim:} {\it In general Im$(\phi)$ is a proper subset of $\Pp\Bb^{id}(M)$.}
\medskip

Here is some examples; to treat them we use some notions introduced in \cite{AGT}, \cite{NA}.

\medskip

(1) There are examples of naked ideal triangulations $T$ of some $\hat
M$ that do not carry any branching at all. It follows from Lemma
\ref{Imphi} that every pre-branching $(T,\omega)$ does not represent a
point in Im$(\phi)$. Let us consider two concrete and instructive
examples.

\medskip

\begin{figure}[ht]
\begin{center}
 \includegraphics[width=9cm]{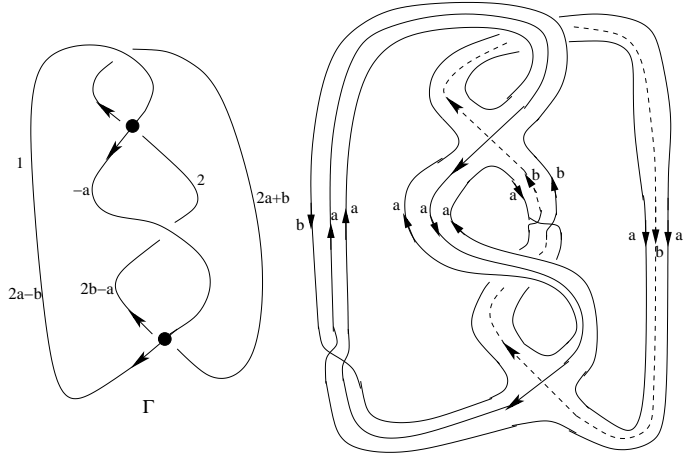}
\caption{\label{8S-1}  No-brancheable, 1 .}
\end{center}
\end{figure}

(a) Let $T$ be the minimal triangulation by two tetrahedra of $\hat M$
where $M$ is the ``compact core'' of the ``figure-8-sister'' cusped
manifold. In Figure \ref{8S-1} we show a weakly branched realization
$(T,\tilde b)$ of it encoded by a $\Nn$-graph $\Gamma$ (see
\cite{AGT}, \cite{NA}).  The ``colors'' $2, 1\in \Z/3\Z \cong A_3$,
$1\leftrightarrow (0,1,2)$, the color $0$ being omitted. The picture
also shows a decoding of $\Gamma$, that is a regular neighbourhood of
Sing$(\Sigma)$ in the dual spine $\Sigma$. We use $\omega_{\tilde b}$
to orient Sing$(\Sigma)$. $\Sigma$ has two $2$-regions on which we fix
auxiliary orientations. The letters $a$ and $b$ refer to the
coefficients of a cellular $2$-chain $\beta$ on $\Sigma$; they label
the two oriented boundaries of the respective $2$-regions. On the
graph $\Gamma$ we indicate the coefficients of $\partial \beta$. Then
one computes immediately that $H_2(\Sigma;\Z/2\Z)\cong
H_2(M;\Z/2\Z)=0$; on the other hand, one realizes that for every
$\Z$-$2$-chain $\beta$ on $\Sigma$, $\partial \beta$ is not
``fundamental'', that is $T$ does not carry any pre-branching $\omega$
such that $[\omega]=0\in H_1(M;\Z)$.

\medskip

\begin{figure}[ht]
\begin{center}
 \includegraphics[width=10cm]{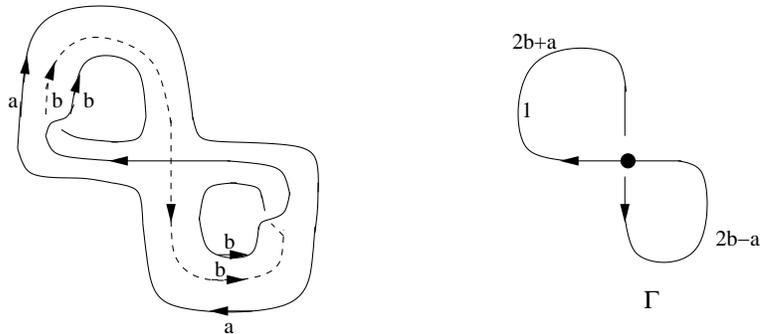}
\caption{\label{nobranching}  No-brancheable, 2.}
\end{center}
\end{figure}

(b) Consider the weakly branched triangulation $(T,\tilde b)$ (of some
$\hat M$) encoded similarly to the above example in Figure
\ref{nobranching}. $T$ consists of one tetrahedron, $\Sigma$ has two
$2$-regions. By easy computations we see that $H_2(M;\Z)=0$,
$H_2(M;\Z/2\Z)\cong \Z/2\Z$. We realize also that every fundamental
$\Z$-$2$-chain $\beta$ on $\Sigma$ is not a branching (that is
$\partial \beta$ is not fundamental). On the other hand, the non
fundamental chain $\beta_0$ such that $a=1$, $b=0$ is a non
orientable branching. Hence, $\omega = \partial \beta_0$ is a
pre-branching on $T$, such that $[\omega]=0\in H_1(M;\Z)$ and
nevertheless $[(T,\omega)]$ does not belong to Im$(\phi)$. One can see
easily that $\partial M$ consists of one spherical boundary component,
so that $\hat M$ is a closed manifold (with a bit of effort one could
realize that $\hat M = \PP^3$ but this is not so important here).

\medskip

(2) Let $(T,\omega)$ be a taut ``layered'' ideal triangulation of
$\hat M$ where $M$ is the compact core of a mapping torus $T_\psi$,
having with fibre a punctured surface $S_V=S\setminus V$ of negative
Euler characteristic. The boundary $\partial M$ is made by tori. Such
a $(T,\omega)$ is constructed by means of an ideal triangulation $K$
of $(S,V)$ (here the compact closed surface $S$ plays the role of
``$\hat S_V$'') and of a sequence of naked flips $K\Rightarrow
\psi(K)$, where $\psi$ is an automorphism of $S$ which pointwise fixes
the punctures $V$ (see \cite{NA} for the details). Then $(T,\omega)$
does not represent a point in Im$(\phi)$. This follows from (1) of
Lemma \ref{h-inv}, because one can check that $[\omega]$ is a non zero
multiple of the dual of the fibre. Fix any $2D$ branching
$(K,\bG)$. By a main result of \cite{2D}, there is a sequence of
$b$-flips $(K,\bG)\Rightarrow (\psi(K),\psi(\bG))$. By using the
corresponding naked sequence we get an instance of layered taut
$(T,\omega)$ as above; by using the $b$-sequence, we get a layered
branched triangulation $(T, b)$. These represent distinct points in
$\Pp\Bb^{id}(M)$.

\medskip

{\bf The effects of the circuit move.}  
Every pre-branching $(T,\omega)$ determines a natural decomposition
of Sing$(\Sigma)$ by {\it oriented circuits}. A {\it circuit move}
on $(T,\omega)$ produces a new pre-branching $(T,\omega')$ by just
inverting the orientation of one circuit. It is clear that every pre-branching
carried by $T$ can be obtained from any given $(T,\omega)$ by performing 
a finite sequence of circuit moves. Hence by adding these moves to the 
ideal $pb$-transits we generate an equivalence relation whose quotient set
consists of one point (see  \cite{BP, AGT}).  We have

\begin{lem}\label{circuit} $\Pp\Bb^{id}(M)$ does not consist of
  one point if and only there are ideal $pb$-triangulations $(T,\omega)$
  and $(T,\omega ')$ of $\hat M$ which differ by a circuit move and
  represent different points in $\Pp\Bb^{id}(M)$.
\end{lem}
\Dim Let $(T_1,\omega_1)$ and $(T_2,\omega_2)$ represent different
points in $\Pp\Bb^{id}(M)$. By Lemma \ref{sameT}, there are  $pb$-transits
$(T_1,\omega_1)\Rightarrow (T,\omega)$, 
$(T_2,\omega_2)\Rightarrow (T,\omega')$.  $(T,\omega)$ and $(T,\omega ')$
represent different points in $\Pp\Bb^{id}(M)$ and are related to each
other by circuit moves, hence at least
one of them changes the class in $\Pp\Bb^{id}(M)$.

\cvd  

\begin{prop}\label{signedQH} Let  $(T,\omega)$ and $(T,\omega')$ be
  triangulations of $\hat M$ which differ by a circuit move supported
  by an oriented circuit $\gamma$ in the $\omega$-oriented
  Sing$(\Sigma)$. Then:

  (i) If $2[\gamma]\neq 0 \in H_1(M;\Z)$, then $(T,\omega)$ and
  $(T,\omega')$ represent different points in $\Pp\Bb^{id}(M)$.

  (ii) Forgetting the orientation, if $[\gamma]_2\neq 0 \in
  H_1(M;\Z/2\Z)$, then $(T,\omega)$ and $(T,\omega')$ represent
  different points in $\Pp\Bb^{id}(M)$.

  (iii) If either $H_1(M;\Z)\neq 0$ and has no $2$-torsion elements,
  or $H_1(M;\Z/2\Z)\neq 0$, then every ideal triangulation $T$ of
  $\hat M$ carries two $(T,\omega)$ and $(T,\omega')$ which differ
  by a circuit move and represent different points in
  $\Pp\Bb^{id}(M)$.
\end{prop}
\Dim (iii) is a consequence of either (i) or (ii). In fact fix any
pre-branching $(T,\omega)$. The $\omega$-oriented circuits of 
Sing$(\Sigma)$ span both
$H_1(M;\Z)$ and $H_1(M;\Z/2\Z)$. At least one, say $\gamma$, is such
that either $2[\gamma]\neq 0$ or $[\gamma]\neq 0$ mod$(2)$.  Implement
the circuit move at $\gamma$ and get the required $(T,\omega')$.

(i) is a direct consequence of (2) of Lemma \ref{h-inv}.

(ii) If $[\gamma]_2\neq 0$, then $[\gamma]\neq 0$. If also
$2[\gamma]\neq 0$, then we are as in point (i). The case when
$2[\gamma]=0$ is less evident. Lift the question in terms of
weak branchings (\cite{NA}, \cite{AGT} , \cite{BP}).
 It is enough to detect a ``character'',
say $\chi$, of every weakly $(T,\tilde b)$ which induces an invariant on
$\Ww\Bb^{id}(M)\sim \Pp\Bb^{id}(M)$ and such that if $(T, \tilde b)$ and 
$(T,\tilde b')$ differ by a circuit move along a circuit that verifies the
conditions in the statement of (ii), then 
they have different $\chi$.  In
\cite{AGT} Remark 8.3 we already noticed that the {\it sign refined}
QH state sums on QH-triangulations with non trivial $c$-weight provide
instances of such characters.

\cvd

The sign refinement mentioned in the proof of (ii) above is
a by-product of \cite{BP}. It would be interesting to produce such a
character $\chi$ in a simpler way without referring to the whole
demanding QH stuff.

 One would wonder that the union of the sufficient conditions
in (i), (ii) of Proposition \ref{signedQH} is also necessary:

\begin{ques}{\rm Let $(T,\omega)$, $(T,\omega')$, $\gamma$ be as
    in the hypotheses of Proposition \ref{signedQH}. Assume that they
    represent different points in $\Pp\Bb^{id}(M)$.  Does it hold true
    that $2[\gamma]\neq 0 \in H_1(M;\Z)$ or $[\gamma]_2\neq 0 \in
    H_1(M;\Z/2\Z)$?}
\end{ques}

Finally let us show that $\Pp\Bb^{id}(*)$ can be infinite.

\begin{prop}\label{infinite-pb} If {\rm rank} $H_1(M;\Z)\neq 0$, then
$\Pp\Bb^{id}(M)$ is infinite.
\end{prop}
\Dim Assume that $\Pp\Bb^{id}(M)$ is finite. Then there is a
triangulation $T$ of $\hat M$ such that every class in
$\Pp\Bb^{id}(M)$ can be represented by some $(T,\omega)$. By varying
the pre-branching $\omega$ we get a finite subset $\Omega =
\{[\omega]\}\subset H_1(M;\Z)$. Take $\gamma \in H_1(M;\Z)$ not
belonging to $\Omega$. There is a simple oriented curve $C$ (not
necessarily connected) traced on the dual spine $\Sigma$ such that
$\gamma = [C]$. By adding along every component $C_i$ of $C$ an
annulus $A_i\subset {\rm Int}(M)$, we get a spine $\SG= \Sigma \cup_i
A_i$ of $M$ which is not even simple. We can thicken the boundary
component of every $A_i$ opposite to $C_i$ to a solid torus to get a
further (non simple) spine of $M$. By applying the ``Bing house''
trick at every solid torus we get a simple spine and possibly using
some ``lune move'' we eventually get a standard spine $\Sigma'$ of $M$
such that the set of homology classes realized by the pre-branchings
of $\Sigma'$ contains $\Omega' = \{\alpha+2\gamma| \ \alpha \in
\Omega\}$. This is absurd.

\cvd

\section{Structures carried by a branched ideal triangulation}
\label{carried}

Let $(M,\partial M)$, $\hat M$ be as usual.
We are going to show that every branched ideal triangulation
$(T,b)$ of $\hat M$ carries a pair of transverse foliations of $M$; 
the $1$-dimensional foliation
will have oriented leaves, and this will be the case also for the $2$-dimensional
provided that $M$ is oriented. 
Every such a foliation can be obtained by integration of some (integrable)
fields of either tangent vectors or (possibly oriented) tangent $2$-planes on $M$. We
will say that two foliations are {\it homotopic} ({\it isotopic}) if
they are obtained by integration of homotopic (isotopic) fields. In
both cases one keeps fields integrability along isotopies. In the case
of vector fields, integrability is preserved also along
homotopies. This is not not required along homotopies of $2$-planes
fields. 

\subsection{The vertical foliation $\Vv_{T,b}$}\label{verticalF} 

\begin{defi} {\rm A {\it traversing foliation} 
$\Ff$ on $M$ is a foliation by {\it oriented} $1$-dimensional 
leaves which satisfies  the following properties:
\begin{enumerate}
\item Every leaf of $\Ff$ is a non degenerale closed interval which
  intersects transversely $\partial M$ at its endpoints.
\item There are {\it exceptional leaves} of $ \Ff$ which
  are simply tangent to $\partial M$ at a finite number of points.
  This tangency points form a compact (non necessarily connected) simple smooth 
  {\it tangency line} $X=X_\Ff$ on $\partial M$. 
  The smooth local model is 
  $$M=\{(x,y,t)\in \R^3\ | t\leq y^2\}, \  \partial M = 
\{(x,y,t)\in \R^3\ | t=y^2\}, \   X=\{t=y=0\}$$ 
  and $\Ff$ is given by the integration of the field 
$\frac{\partial}{\partial t}$ restricted to $M$. 
  \item  $\Ff$ is {\it generic} if every {\it generic} exceptional 
leaf is by definition tangent
  to the boundary at one point, and possibly there is a finite set 
  of {\it non generic} exceptional leaves which are tangent at $2$ points
  {\it with transverse tangency lines}: this means that they are transverse
  in $\partial M$ provided that a neighbourd of one point is transported
  onto a neighbourhood of the other point by using the flow of the traversing
  foliation.
  \end{enumerate}
}
\end{defi}

Since \cite{BP2} we have pointed out (for oriented $M$ and mainly in
terms of {\it oriented} branched spines) that every ideal branched
triangulation $(T,b)$ of $\hat M$ carriers a generic traversing
foliation of $M$ called here the {\it vertical foliations} $\Vv =
\Vv_{T,b}$. In fact this holds as well if $M$ is not orientable.  We
outlines its systematic construction in the form of a {\it puzzle}.
\smallskip

{\bf $\Vv$-puzzle.} Let $T$ be an ideal triangulation of $\hat M$.
Associated to the family of abstract tetrahedra $\{ \Delta\}$ which
form $T$ there is a family of {\it truncated tetrahedra} $\{\tilde D
\}$. The boundary of every such a $\tilde \Delta$ has four triangular
faces (corresponding to the four truncated vertices of $\Delta$) and
four hexagonal faces (the truncation of the four $2$-faces of
$\Delta$). The gluing in pairs of the abstract $2$-faces of $\{
\Delta\}$ which produces $T$ restricts to a gluing in pairs of the
hexagonal faces of $\{\tilde D \}$ which produces a cell decomposition
$(\tilde T, \partial \tilde T)$ of $(M,\partial M)$. The restriction
$\partial \tilde T$ is in fact a triangulation of $\partial M$ made by
the triangular faces of $\{\tilde D \}$.

Recall now a classical way to prove that the Euler-Poincar\'e
characteristic $\chi(Y)$ of a boundaryless compact manifold $Y$
defined by means of the zero indices of any tangent vector fields on
$Y$ with isolated zeros coincides with its combinatorial definition in
terms of any triangulation of $Y$ (having only material vertices):
given such a triangulation $T$ we take its first barycentric
subdivision $T^{(1)}$ endowed with a standard $\Delta$-complex
structure so that in particular every vertex of $T$ is a pit. Then every
symplex of $T^{(1)}$ (with ordered vertices) carries a so called {\it
  Whitney tangent vector field} which can be defined explicitely in
terms of its barycentric coordinates (see \cite{HT}).  All these
locally defined vector fields match to define a globally defined
vector field on $Y$ with an isolated zero at each barycenter of the
`geometric' symplexes of $T$ (with not necessarily ordered vertices)
in such a way that the sum of these indices equals the combinatorial
characteristic of $T$.

\begin{figure}[ht]
\begin{center}
 \includegraphics[width=8cm]{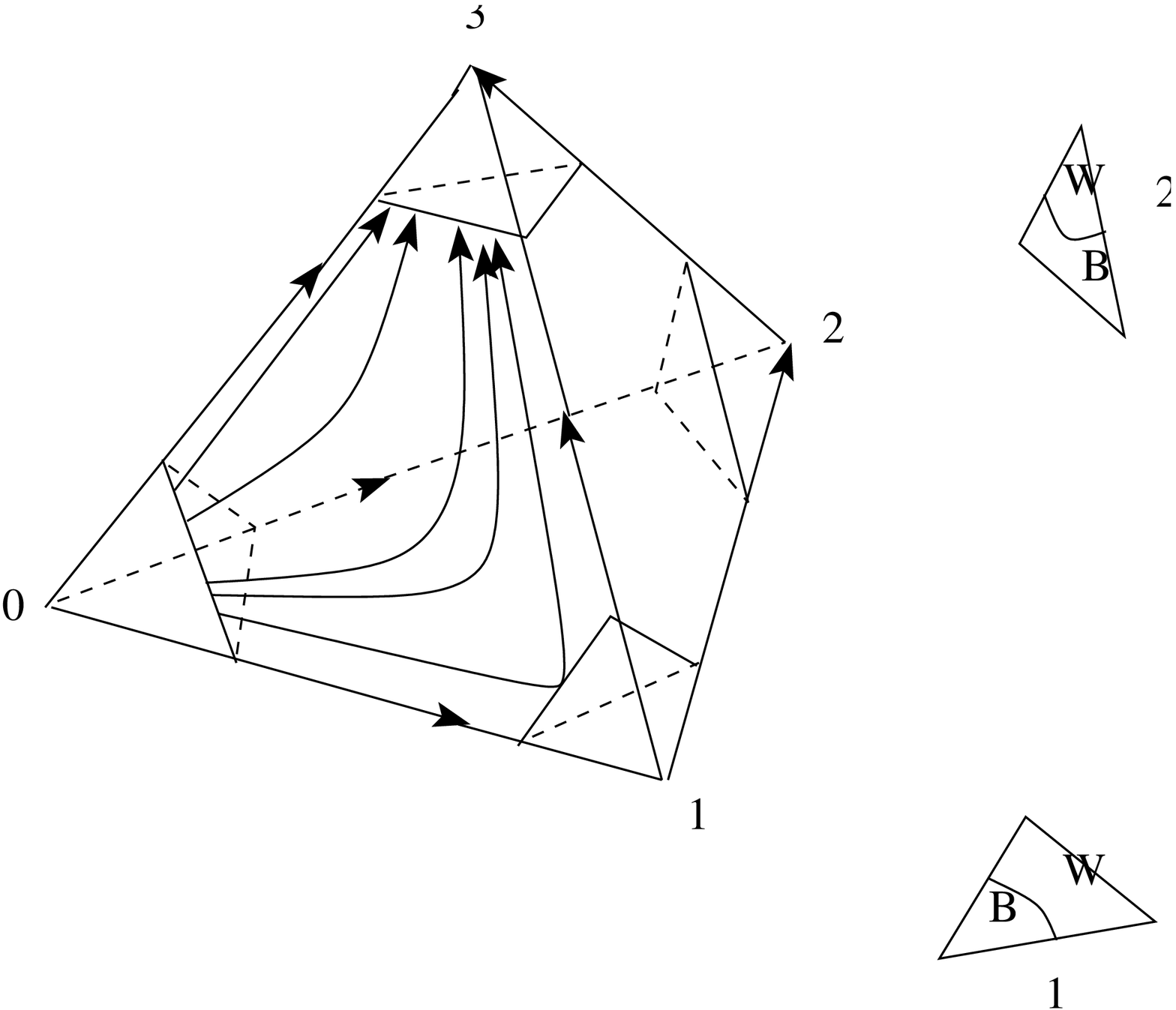}
\caption{\label{Vpuzzle2} Tiles for $(\Vv,\partial \Vv)$.}
\end{center}
\end{figure}

From this general construction we retain that every branched
tetrahedron $(\Delta, b)$ carries such a Whitney field which restricts
to its $2$-dimensional instance on every branched $2$-face of
$(\Delta,b)$ (see Figure \ref{Vpuzzle2} and also the top of Figure
\ref{2Dpuzzle}).  Finally we consider the restriction of the oriented
integral lines of these fields to every truncated tetrahedron $\tilde
\Delta$.  This gives us the tiles $(\tilde \Delta, \Vv_{\tilde
  \Delta})$ of our puzzle. For every branched triangulation $(T,b)$
of $\hat M$, they match to produces the required vertical traversing
foliation $\Vv_{T,b}$.

\smallskip

{\bf $\Vv$-boundary bicoloring}. Every traversing foliation $\Ff$ of
$M$, in particular a vertical foliation $\Vv=\Vv_{(T,b)}$, determines
a {\it bicoloring}, denoted by $\partial \Ff$, of the components of
$\partial M \setminus X_\Ff$: let us say that a component $C$ is {\it
  white} ({\it black}) if the foliation is ingoing (outgoing) at
$C$. For every tile $(\tilde \Delta, \Vv_{\tilde \Delta})$ denote by
$t_j$, $j=0,1,2,3$, its triangular face at the truncated vertex $v_j$
of $(\Delta, b)$ (as usual the vertices of $(\Delta,b)$ are ordered
via a labelling by $0,1,2,3$); then $t_0$ is all white, $t_3$ is all
black; both $t_1$ and $t_2$ are divided by an arc of $X_\Vv$ into a
white and a black portion; the white (black) zone of $t_1$ ($t_2$)
contains two vertices of that triangular face (see again Figure
\ref{Vpuzzle2}).  Denote by $W_b(\partial M)$ (resp. $B_b(\partial
M)$) the union of the white (black) regions of $\partial M$ determined
by the boundary bicoloring $\partial \Vv$.  Notice that
$$\chi(M)=  \frac{1}{2}\chi(\partial M) = \chi (\overline {W_b(\partial M)})
= \chi (\overline {B_b(\partial M)})\ . $$  

{\bf $M$ oriented.}  If $M$ is {\it oriented} the colors can be encoded by an
orientation: a black component keeps the
boundary orientation (according to the usual rule ``first the
outgoing normal''), while a white component has the opposite
orientation.

When $M$ is oriented there is another systematic way to construct
$\Vv$. Let $M$, $(T,b)$, $(\Sigma, b)$ be
as usual. Let
$(\Sigma^*,b^*)$ be an ``abstract'' copy of the oriented branched
surface $(\Sigma,b)$, with its own branched smooth structure. Let us
consider the oriented branched $3$-manifold with
boundary $$F=F(T,b)=(\Sigma^* \times I, \hat b), \ I=[-1,1] \ . $$ Its
singular set is
$$ {\rm Sing}(F)= {\rm Sing}(\Sigma^*)\times I$$ hence $F$ has no
``vertices'', i.e. $0$-dimensional singular strata.  For every edge
$e$ of Sing$(\Sigma^*)$, $\{e\}\times I$ is called a {\it switching
  surface} of $F$. For every vertex $v$ of Sing$(\Sigma^*)$,
$\{v\}\times I$ is called a {\it pivot} of $F$. The orientation $\hat
b$ on every $3$-dimensional region of $F\setminus {\rm Sing}(F)$ is
the product of the natural orientation of $I=[-1,1]$ by the
$b^*$-orientation of a $2$-region of $\Sigma^*$. Note that we have
inverted the factors order. Switching surfaces and pivots are oriented
similarly.  $F$ carries the ``vertical'' foliation $\Vv^*$ by the
oriented segments $\{x\} \times I$.  Every pivot, swithching surface,
or $3$-dimensional region of $F$ is union of leaves of $\Vv^*$. An
interesting result of \cite{GR} is that $M$ can be (piecewise
linearly) embedded into $F$. By taking a suitable {\it normal smooth}
embedding we can realize $\Vv_{T,b}$ as the restriction of $\Vv^*$ to
$M$. This idea had been already exploited in \cite{BP3}; we will add
more precision about a systematic construction of such a normal
embedding, again in the form of a puzzle. Before doing it, we
summarize its main features.

\begin{prop}\label{M-emb-vert} For every branched triangulation $(T,b)$
there is a normal smooth embedding  $M\subset F(T,b)$
which satisfies the following properties:
\begin{enumerate}

\item The restriction $\Vv= \Vv_{T,b}$ of the vertical foliation
  $\Vv^*$ of $F=F(T,b)$ is a generic traversing foliation of $M$ which
  coincides with the one described via the $\Vv$-puzzle.

\item Denote by $M^\circ$ the maximal subset of $M$ formed by
vertical segments which are leaves of both $\Vv^*$ 
and $\Vv$. Then $M^\circ$ includes Sing$(F)$ and $F\setminus (\Nn \times I)$
where $\Nn$ is a regular neighbourhood $\Nn$ of Sing$(\Sigma^*)$
in $\Sigma^*$.

\item The vertical leaves contained in Sing$(F)$ coincide with the
exceptional leaves of $\Vv$.  The pivots of $F$ are  
the non generic exceptional leaves. 

\item The oriented branched surface $\Sigma \subset M \subset F$ is
  transverse to the leaves of $ \Vv$, in such a way that:
  
  (i) It intersects all leaves.
  
  (ii) Every generic exceptional leaf intersects $\Sigma$ at $2$ points;
  every pivot intersects $\Sigma$ at $4$ points.
  
  (iii) Every leaf which intersects Sing$(\Sigma)$ is non exceptional
   and intersects $\Sigma$ at one point. 
   
  (iv)  The product
  orientation $I\times b$ coincides with $\hat b$.

\item The normal embeddings of $M$ in $F$ are considered up to 
isotopy through normal embeddings, so that $\Vv$ is uniquely
determined up to isotopy, preserving the exceptional leaves. 
\end{enumerate}
\end{prop}

\cvd

\begin{figure}[ht]
\begin{center}
 \includegraphics[width=10cm]{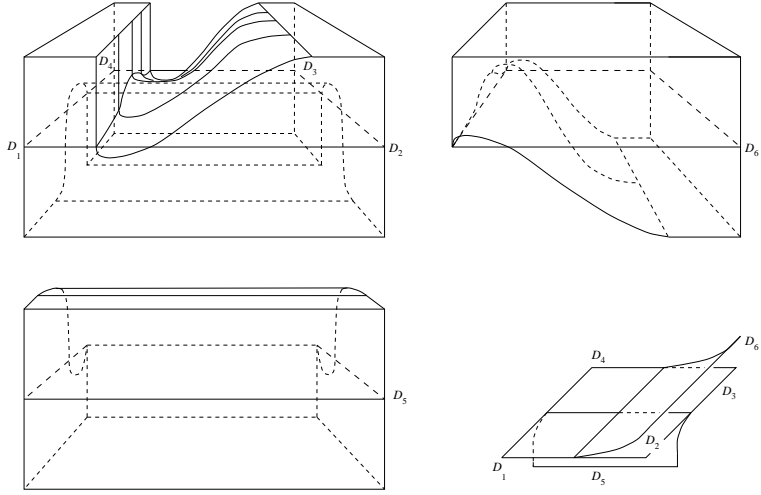}
\caption{\label{figure11} Tiles for the normal embedding of $M$ into $F(T,b)$.}
\end{center}
\end{figure}

{\bf The normal embedding of $M$ into $F(T,b)$}
The construction below is illustrated in Figure \ref{figure11}.

(1) We decompose $\Sigma^*$ (hence $\Sigma$) by ``tiles'' of three
types: {\it butterfly}, ``$Y\times J$'', and {\it disk}.  The
butterflies form a branched regular neighbourhood of the vertices of
$\Sigma^*$. Every butterfly has a boundary made by the union of four
branched tripods and six simple arcs glued at eleven corners at which
an arc and a tripod have a common endpoint.  Every $Y\times J$ is the
product of a branched tripod $Y$ and a closed interval $J$. Its
boundary is the union of two tripods and three simple arcs, with six
corners. There is one such a tile at every edge of
Sing$(\Sigma^*)$. Butterflies and $Y\times J$ tiles glue at pairs of
common boundary tripods to form a branched regular neighbourhood $\Nn$
of Sing$(\Sigma^*)$ in $\Sigma^*$.  The boundary $\partial \Nn$ is a
union of smooth simple circles. The closure of every component of
$\Sigma^* \setminus \Nn$ is a disk tile. The union of these disks
intersect $\Nn$ along $\partial \Nn$.

\medskip

(2) Let us generically denote by $P$ such $2$-dimensional tiles.  Then
$P\times I$ is the corresponding tile of $F=\Sigma^* \times I$.  These
$3$-dimensional pieces match at the surfaces $\partial P \times I$ so
that their union reconstructs the whole of $F$. Within every $P\times
I$ we will specify a smoothly embedded $3$-dimensional manifold $Q_P $
(with boundary and corners), such that $Q_P$ also embedds into $M$ and
the above matching coherently restricts to the surfaces $Q_P \cap
(\partial P \times I)$; the union of the $Q_P$'s eventually
realizes the required embedding of $M$ into $F$. We stress that every
$Q_P$ will be determined by the actual embedding of (a copy of) $P$
into $M$, by using both the branching $b$ and the orientation of
$M$. Let us describe $Q_P$, type by type.

$\bullet$ Let $P$ be a disk. Then $Q_P$ coincides with the whole of
$P\times I$.  The embedding of $Q_P$ in $M$ is such that $P\times
\{0,\}$ is the copy of $P$ in $\Sigma$.

\medskip

$\bullet$ Let $P$ be a tile of type $Y\times J$. Then also $Q_P$ is a
product of the form $S_P\times J$, where $S_P \subset Y\times I$, and
for every $s\in J$, $S_P\times \{s\} \subset Y\times \{s\}\times I$
fiberwise.  Let us denote by $l_-,\ l_+$ the two legs of the oriented
branched tripod $Y$ which induce the prevalent orientation on the
central vertex $x_0$.  The order of the legs is determined by the
$b$-orientation of of Sing$(\Sigma^*)$ and the orientation of $M$. Let
$l_0$ be the other leg.  Let us subdivide $l_\pm$, by its midpoint
$x_\pm$. Hence $l_\pm=l_\pm'\cup l"_\pm$, where $l_\pm'=[z_\pm,x_\pm]$
and $l"_\pm=[x_\pm,x_0]$ Set $y_\pm$ the midpoint of $l"_\pm$, $w$ the
one of $l_0=[x_0,w]\cup [w,p]$.  On $l_\pm \times I$ consider function
$t=f_\pm(x)$ such that: 

\medskip

$f_\pm$ is smooth on $l_\pm \setminus \{x_0\}$;

\medskip

$ f_\pm(x)=\mp1$,  if $ x\leq x_\pm$;

\medskip

$f_+$ (resp. $f_-$) is incresing (decreasing) on $[x_\pm,y_\pm]$,
decreasing (increasing) on $[y_\pm,x_0]$;

\medskip

$f_\pm(y_\pm)=\pm 1/4$;

\medskip

$f_\pm(x_0)=0$, and the graph of $f_\pm$ is simply tangent to the vertical
segment $\{x_0\}\times I$ at the point $(x_0,0)$.

\bigskip

Set 
$$S_P^+ = \{(x,t)| \ t\geq f_+(x)\}, \ S_P^- = \{(x,t)| \ t\leq f_-(x)\} \ . $$

Finally set

$$S_P= S_P^+ \cup S_P^- \cup (l_0\times I) \ . $$
Note that the boundary of $S_P$ contains a smooth line, $\gamma$ say,
formed by the union of the graphs of the functions $f_\pm$ which is
simply tangent to $\{x_0\}\times I$; hence the surface $\gamma \times
J$ is simply tangent to $\{x_0\}\times J \times I$ along
$\{x_0\}\times J\times \{0\}$. The embedding of $Q_P$ into $M$ is such
that the copy of $P$ in $\Sigma$ is subdivided by: \medskip

$P\cap \{(l_\pm'\times J)\times  I\}= (l_\pm' \times J)\times \{0\}$;

\medskip

$P\cap \{([w,p]\times J)\times I\}= ([w,p] \times J)\times \{0\}$;

 \medskip
 
 $P\cap \{((S_P^\pm\cup [x_0,w])\times J)\times I\}$ is ``parallel''
 to $\{\gamma \cup ([x_0,w])\times \{0\}) \} \times J $ with switching
 curve given by $\{w\}\times J \times \{0\}$.
 
 \bigskip

 $\bullet$ Let $P$ be a butterfly. We can normalize the picture so
 that $P$ is formed by a {\it plate} $\Pp$ and two {\it wings}
 $\Ww_\pm$.  The wings are ordered by the $b$-orientation of $\Pp$ and
 the orientation of $M$; within $M$ it makes sense that $\Ww_+$ lies
 ``over'' $\Ww_-$.  We can also assume that $\Pp$ is obtained by
 removing from a square $Q=J^2$ (with coordinate $(x,y)$) four
 disjoint sectors of open $2$-disks centred at the corners of $Q$; the
 boundary of every sector intersects $\partial Q$ at an arc properly
 embedded in $Q$ which near $\partial Q$ is made by two orhogonal
 small segments.  The four tripods in $\partial P$ have centers in the
 midpoints of the edges of $Q$, two legs contained in $\partial \Pp$
 and a further leg contained in $\partial \Ww_\pm$. The union of the
 wings intersects $\Pp$ at the cross $\{xy=0\}$ (which is included in
 Sing$(\Sigma^*)$). The $b$-orientation of the cross determines two
 bands $\Bb_\pm$ of $\Pp$ ``on the left'' of either $\{x=0\}$ or
 $\{y=0\}$.  These are ordered according to the fact that the wing
 $\Ww_\pm$ folds over $\Bb_\pm$ in $\Sigma \subset M$.  To fix the
 idea, let us assume that $P$ corresponds by duality to a tetrahedron
 of $(T,b)$ such that the $b$-sign $*_b=+1$ (the other case can be
 treated similarly). Then $\Bb_+$ is bounded by $\{y=0\}$, while
 $\Bb_-$ is bounded by $\{x=0\}$.  Take $\Pp\times I \subset F$. This
 also embedds into $M$. By adopting the above notations, we can ``dig
 a groove'' in $\Bb_+ \times I$ and realize an embedding of $(S_P^-
 \cup (l_0\times I))\times J$ into $\Pp\times I$ (hence in $M$), so
 that $\{y=0\}\subset \Pp$ is the tangency line.  Then we can glue a
 copy of $S_P^+ \times J$ along $\{y=0\}\times [0,1]$. We call this a
 {\it thick wing}.  The resulting space $Q_P^+$ embedds in both $F$
 and $M$. $Q_P^+$ can be glued to suitable pieces $Q_*$ constructed as
 above along two edges (possibly the same) of Sing$(\Sigma^*)$. We can
 manage similarly (and independently) on $\Bb_-$, dig a groove to
 realize an embedding of $(S_P^+ \cup (l_0\times I))\times J$ and
 complete it by a thick wing $S_P^- \times J$ along $\{x=0\}\times
 [-1,0]$. This produces $Q_P^-$. For every $\Pp$ consider the {\it
   intersection} $Q_P'=Q_P^+\cap Q_P^-$. These can be assembled with
 the $3$-dimensional pieces of the other types to get a spaces $M'$
 which embedds in both $F$ and $M$; $M\setminus M'$ consists of the
 union of disjoint small neighbourhoods of the vertices of $\Sigma$
 (each one corresponding to the portion of $\Pp \times I$ where the
 two grooves cross).  In order to get the whole of $M$, we define the
 ultimate $Q_P$ by modifing by isotopy the two grooves (and
 consequently the two thick wings) in such a way that $Q_P$ coincides
 with $Q_P'$ near the boundaries and the grooves do not
 cross. Precisely, the new groove digged in $\Bb_+ \times I$ has
 tangengy line on $\{y=0\}\times [0,1]$ which is the graph of a smooth
 non negative bell function $t=g_+(x)$ which is equal to $0$ near the
 boundary, and $g_+(0)=2/4$ is its maximum value. Similarly, the new
 groove digged in $\Bb_- \times I$ has tangengy line on $\{x=0\}\times
 [-1,0]$ which is the graph of a smooth non positive bell function
 $t=g_-(y)$ which is equal to $0$ near the boundary, and $g_-(0)=-2/4$
 is its minimum value. The grooves are isotopically either lifted up
 or lowered consequently.  The embedding of a butterfly of $\Sigma$
 into such a $Q_P$ agrees with the embedding described above near the
 boundary and can be naturally extended in the interior.

\medskip

We have achieved the promised {\it normal smooth embedding} of $M$
into $F(T,b)$.  By construction a normal embedding preserves the
orientation. By construction we know exactly the intersection of
$M^\circ$ with the pieces $Q_P$ of any type.

\subsection{The horizontal foliation $\Hh_{T,b}$}\label{horizontalF} 
Let $X=X_{\Vv}\subset \partial M$ be the system of tangency lines
of a  vertical foliation $\Vv=\Vv_{T,b}$  as above. 
We can thicken every component $C$ of $X$ to an annulus $A_C$
in $\partial M$, foliated by parallel copies of $C$. This gives us a system
of {\it sutures} on the boundary $\partial M$. Denote by $\Aa$ the union of
these annuli. The horizontal foliation
$\Hh=\Hh_{T,b}$ has the following main properties:

\begin{enumerate}

\item $\Vv$ and $\Hh$ are transverse foliations.

 \item The
  closure of every component of $\partial M \setminus \Aa$ is a leaf
  of $\Hh$, while $\Hh$ is
  transverse to $\partial M$ along $\Aa$ and induces on every 
  $A_C$ the prescribed foliation. This is called the boundary configuration
  $\partial \Hh$ of $\Hh$.
  
  \item $\Hh$ is uniquely determined up to homotopy keeping the boundary
  configurations up to isotopy.
\end{enumerate}
  
If $M$ is oriented, $\Hh$ is oriented as well in such a way that $\Vv$
intersects it everywhere with intersection number equal to $1$.
Also $\Hh$ can be produced by a puzzle.

\begin{figure}[ht]
\begin{center}
 \includegraphics[width=10cm]{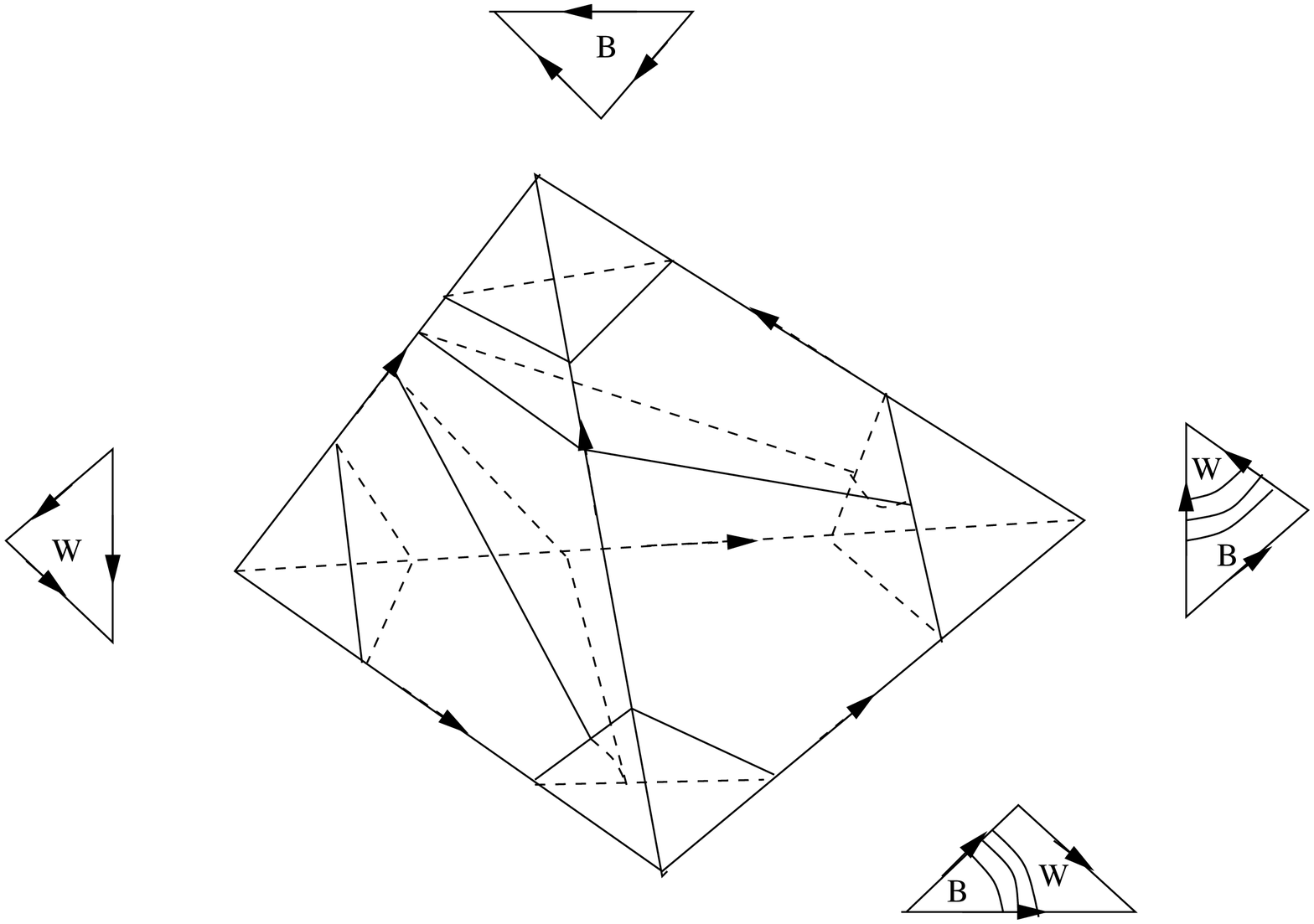}
\caption{\label{3Dpuzzle} $3D$ tiles for $\Hh$.}
\end{center}
\end{figure}

\begin{figure}[ht]
\begin{center}
 \includegraphics[width=4cm]{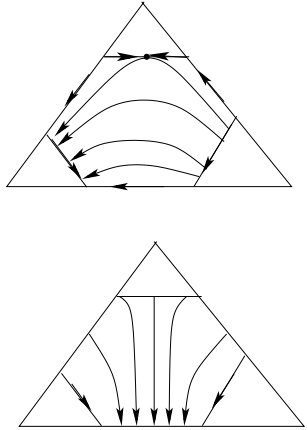}
\caption{\label{2Dpuzzle} $2D$  tiles.} 
\end{center}
\end{figure}

{\bf $\Hh$-puzzle.} Every truncated tetrahedron $\tilde \Delta$
associated to $(\Delta, b)$ carries a tile $(\tilde \Delta,
\Hh_{\tilde \Delta})$.  This is illustrated in Figure
\ref{3Dpuzzle}. The picture shows also the boundary configurations
which is in agreement with the tiles of $\partial \Vv$ considered
above (including the 2D branching $\partial \omega_b$ provided that
$M$ is oriented - let us forget it in general). In the triangles $t_1$
and $t_2$ we see the trace of $\Aa$ and the colored complementary
portions. Every colored boundary piece is part of some horizontal
leaf. In the bulk we see two typical leaves of $\Hh_{\tilde \Delta}$
which are transverse to $\partial M$ along $\Aa$. The trace of
$\Hh_{\tilde \Delta}$ on each hexagonal face of $\tilde \Delta$ is a
tile of the $2D$ analogous of $\Hh_{\tilde \Delta}$ (see \cite{2D}).
This is illustrated on the bottom of Figure \ref{2Dpuzzle}.
  
If $M$ is oriented, also $\Hh=\Hh_{T,b}$ can be produced starting from
a normal smooth embedding of $M$ into $F=F(T,b)$. $F$ carries also a
horizontal foliation $\Hh^*$ with oriented branched leaves of the form
$\Sigma^* \times \{t\}$. $\Hh$ is not immediately equal to the
restriction $\Hh'$ of $\Hh^*$ to $M$. They coincide on $M^\circ$, but
$\Hh'$ has a somewhat complicated behaviour at $\partial M$ on
$M\setminus M^\circ$. Eventually $\Hh$ is homotopic to $\Hh'$.

\section{On $\Ss^{id}(*)$}\label{slide}
Now we can understand the meaning of the combinatorial 
classification of ideal $b$-transit (see Section \ref{sub-fam}).
We mainly refer to \cite{BP4} \cite{BP2}.
We have
\begin{prop}\label{sliding} 
Let $(T,b)\to (T',b')$ be an ideal ($2 \leftrightarrow 3$
or quadrilateral $0\leftrightarrow 2$) branched transit. 
Then the following facts are equivalent to each other:

\medskip

(1) It is a sliding transit;
\medskip

(2) $\partial \Vv_{(T,b)}= \partial \Vv_{(T',b')}$, that is the
boundary bicoloring of the vertical foliations is preserved
(up to isotopy);
\medskip

(3) There is a smooth isotopy $\Psi: M\times [0,1]\to M$, 
$\psi_t= \Psi_{|M\times \{t\}}$, such that:
\begin{enumerate}
\item[(i)]  $\psi_0=id$; for every $t>0$, $\psi_t$ is an embedding
onto $M_t \subset {\rm Int}(M)$ and $M\setminus M_t$ is a collar of
$\partial M$.
\item[(ii)] For every $t$, the restriction of $\Vv_{(T,b)}$ to $M_t$
is a traversing foliation $\Ff_t$ which is not generic at one 
value  $t_0\in (0,1)$; hence the pullback $\Vv_t:= \phi_t^*(\Ff_t)$ is a traversing
foliation on $M$.
\item[(iii)] For every $t\in [0,t_0)$, $\Vv_t$ is isotopic to $\Vv_{(T,b)}$;
for every $t\in (t_0,1]$, $\Vv_t$ is isotopic to $\Vv_{(T',b')}$ 
\end{enumerate} 
In particular $\Vv_t$ is a homotopy through
traversing foliations that connects $\Vv_{(T,b)}=\Vv_0$ with $\Vv_{(T',b')}$.
\end{prop}

Then implications ``$(3)\Rightarrow (2)$'' holds in general for every
homotopy through traversing foliations: such a homotopy induces an
isotopy of the boundary bicolorings (see \cite{BP2}, Lemma 4.3.8).  In
such a generality ``$(2)\Rightarrow (3)$'' fails, but it is true in
the restrictive hypothesis of the Proposition. This can be checked
move by move; see again \cite{BP2}.
\smallskip

 Point (a) of the following Proposition is equivalent to Proposition
 \ref{sliding} (2); point (b) adds information about the way the
 boundary bicoloring changes.  For a proof see \cite{BP2}, Proposition
 3.5.1.
  
\begin{prop}\label{bump}  Let $(T,b)\to (T',b')$ be an ideal  branched transit. Then:

(a) the following facts are equivalent to each other:
\begin{enumerate}

\item It is a bump transit;
\medskip

\item $\partial \Vv_{(T,b)}\neq \partial \Vv_{(T',b')}$
\end{enumerate}

(b) Precisely the bicolorings differ to each other by an instance of
the following modifications:
\begin{enumerate}
\item[(i)] $\overline{B_{b'}(\partial M)}$ is obtained by adding to
  $\overline{B_{b}(\partial M)}$ one $1$-handle embedded in
  $\overline{W_{b}(\partial M)}$ and one $2$-disk embedded in
  $W_{b}(\partial M)$, disjoint from the attached handle;

\item[(ii)] The inverse of the modification in (i);

\item[(iii)] The modifications obtained as above provided that the roles of the white/black
portions of $\partial M$ are  exchanged.
\end{enumerate}
\end{prop}

\begin{remark}{\rm  The above bicoloring modifications preserve the necessary property that
$$ \chi (\overline {W_b(\partial M)})
= \chi (\overline {B_b(\partial M)})\ . $$}
\end{remark}

We are ready to point out the intrinsic content of
$\Ss^{id}(M)$.
 
\begin{prop}\label{S-char} Let $(T,b)$ and $(T',b')$ be ideal branched
  triangulations of $\hat M$. Then:
  
  The following facts are equivalent to each other:
  
  (1) They represent the same point in
  $\Ss^{id}(M)$;
  \medskip
  
  (2)  $\Vv=\Vv_{(T,b)}$ and
  $\Vv'=\Vv_{(T',b')}$ are homotopic through traversing
  foliations.
  \medskip
  \end{prop}

\smallskip

``$(1)\Rightarrow (3)$'' by Proposition \ref{sliding}.
Clearly  ``$(3)\Rightarrow (2)$''. The proof of ``$(2)\Rightarrow (3)$'' is
essentially equivalent to the proof of  Theorem 4.3.3.  of \cite{BP2}, 
see especially Proposition 4.4.9. There one considers closed oriented $M$ with 
triangulations with only one vertex, but the proof runs in general without substantial differences. 
It deals in terms of the dual branched spines
which at present are transversely oriented, not necessarily oriented.
Actually the proof in easier here
because we are allowing arbitrary sliding lune 
moves, so that we can avoid most of the discussion of Section 4.5. of \cite{BP2}.
Let us recall a few points of the proof: by trasversality, we can assume 
that the homotopy is  {\it generic}; this means that every $\Vv_t$ is generic 
with the exception of a finite number of points $t_j\in (0,1)$,
 $j=0,\dots,n$ where the genericity can be lost according to a determined
 finite set of configurations.  The boundary bicoloring is constant along the 
 whole homotopy; for every generic $\Vv_t$ the traversing foliation is  carried
 by  transversely  oriented branched {\it simple} spine $\Sigma_t$, such that   
 $\Sigma_0=\Sigma$, $\Sigma_1=\Sigma'$ (hence they are standard); 
 along every interval $(t_j,t_{j+1})$,  the foliations as well as the spines are 
 isotopic to each other, and  by analyzing the possible configurations, 
 one realizes that passing through every special value $t_j$ has the 
 effect to perform locally a sliding move (which makes sense also for
 simple, not necessarily standard, spines). It might actually happen that
 the standard setting is lost at some event where a negative lune move is 
 performed; however by using other sliding moves (modifying the homotopy
 itself) we can overcome it restoring the standard setting everywhere. 
 
 \cvd
 
 \medskip
 
 Let us say that a homotopy as in (3) of Proposition \ref{S-char}
 {\it realizes the sliding equivalence}.
 
 \begin{cor} If  $(T,b)$ and $(T',b')$ represent the same point
 in $\Ss^{id}(M)$, a homotopy that realizes the sliding equivalence
 can be augmented to a homotopy of couples of trasverse foliations   
 $(\Vv_t, \Hh_t)$  which induces an isotopy of the couples of boundary 
 structures $(\partial \Vv_t, \partial \Hh_t)$, that is  it connects
 the patterns of structures carried by $(T,b)$ and $(T',b')$ respectively. 
 \end{cor}
 \smallskip
 
 We can strengthen the intrinsic content of the sliding equivalence by considering
more general $1$-foliations on $M$. Let $\Ff$ be a non singular foliation
 on $M$ by oriented curves which along the boundary  $\partial M$ has the 
 same qualitative features of a traversing foliation but it is not necessarily a
 traversing foliation. 
 Denote by $\partial \Ff$ its boundary
 configuration, that is the usual bicoloring. Call {\it admissible}
 every boundary configuration $\fG$ obtained in this way.
 We have (see \cite{BP4}):
 \begin{lem} A boundary configuration $\fG$ is admissible
 if and only if $\chi(M)=\chi(\overline{W(\fG}))=  \chi(\overline{B(\fG}))$.
 \end{lem} 
 
 Fix such an admissible boundary configuration $\fG$.
 Denote by $\Ff(M,\fG)$ the set of non singular foliations 
 $\Ff$ on $M$ such that $\partial \Ff = \fG$ (up to isotopy), 
 considered up to homotopy through non singular foliations which 
 is an isotopy at $\partial M$.  We stress that we are not requiring 
 that $\Ff$ is traversing. 
 
 Let $M'$ be obtained by removing an open $3$-ball from $M$,
 creating a new spherical boundary component $S$.
 Given $\fG$ as above, denote by $\fG'$ the admissible 
 boundary configuration at $\partial M'$ that extends
 $\fG$ in such a way that 
 the bicoloring on $S$ consists of one black and one white disk
 separated by one simple curve.
 Denote by
 $\Ss^{id}(M',\fG')$ the subset of $\Ss^{id}(M')$ formed by
 the classes of triangulations $(T,b)$ of $\hat M'$ such that 
 $\partial \Vv_{(T,b)} = \fG'$. Then,
 up to homotopy, the traversing foliation $\Vv_{T,b}$ is the restriction of
 some non singular foliation representing an element of  
 $ \Ff(M,\fG)$. Hence there is a well defined map
$$\psi:   \Ss^{id}(M',\fG') \to \Ff(M,\fG) \ . $$
 A main result of \cite{BP4} can be rephrased as follows.

 \begin{prop}\label{flow-spine} For every admissible boundary
configuration $\fG$ at $\partial M$, the map
$$\psi:   \Ss^{id}(M',\fG') \to \Ff(M,\fG) $$
is bijective. 
\end{prop}

The particular case when $M'$ is oriented and $\partial M'$ consists of
just one spherical boundary component $S$ had been early considered in
\cite{BP2}.  In this case $M$ is a closed manifold and
$\Ff(M,\emptyset)$ is the set of arbitrary non singular
foliations of $M$ by oriented curves, considered up to homotopy
sometimes called ``combings''; the generalization in \cite{BP4} 
(again for oriented $M$) is not hard.
He we consider also $M$ non orientable but the proof holds as well. 
First one proves at the same time that  $\Ss^{id}(M',\fG')$
is non empty and the map $\psi$ is onto (see Proposition
5.1.1 of \cite{BP}). The proof is based on Ishii's notion of {\it flow
  spines} \cite{I}. For the injectivity of the map, see Theorem 5.2.1 of \cite{BP}. 
  The basic idea is to `cover' any homotopy with a
chain of flow-spines connecting $(T,b)$ with $(T',b')$ such that the
traversing foliation associated to one is homotopic through traversing
foliations to the traversing foliations of the subsequent.

\subsection{An alternative approach to the  connectivity results}\label{alternative}
By using the above results, the natural projection $\Ss^{id}(M)
\to \Bb^{id}(M)$ and the fact that $\Bb^{id}(M)$ consists of one point
by Theorem \ref{bI}, can be rephrased by saying that the bump moves
modify in a transitive way the traversing foliations carried by the
ideal branched triangulations of $\hat M$. By elaborating on this
remark, we will outline a (partially conjectural) different, perhaps more conceptual
approach to a proof of the main branched connectivity
result and ultimately also of the naked one. With the notations of
Proposition \ref{flow-spine}, let us fix an admissible
boundary configurations $\fG$ and $\fG'$ of $M$ and $M'$. We have

\begin{prop}\label{pont}
  Let $(T,b)$ and $(T',b')$ be branched ideal triangulations of $\hat M'$
  representing two different points of $\Ss^{id}(M', \fG')$. Then
  there is a composite ideal $b$-transit $(T,b)\Rightarrow (T',b')$
  (necessarily including bump moves).
\end{prop}
\smallskip

Obviously this follows from Theorem \ref{bI} but we can prove it
independently.  $\Ss^{id}(M',\fG') \sim \Ff(M,\fG)$ is an affine space
on $H_1(M; \Z)$ (hence they are in general infinite sets). 
Two non singular foliations on $M$ with the same
boundary configuration differs to each other (up to homotopy through
foliations of the same type) by a so called {\it Pontrjagin move};
this can be realized by a so called {\it combinatorial Pontrjagin
move} on $(T,b)$ which is realized by a composite ideal $b$-transit as in the
statement of the Proposition; details can be found in the proof of
Proposition 6.3.1 of \cite{BP2}.

\begin{prop}\label{move-col} Let $(T_1,b_1)$ be an ideal branched triangulations  
  of $\hat M'$ representing a point of $\Ss^{id}(M', \fG_1')$
  and $\fG_2$ be another admissible boundary configuration for $M$.
  Then there is a composite ideal $b$-transit $(T_1,b_1)\Rightarrow (T_2,b_2)$
  such that $(T_2,b_2)$ represents a point in $\Ss^{id}(M', \fG_2')$.
\end{prop}

\smallskip

Again it is a consequence of Theorem \ref{bI}. However the statement is apparently
weaker. Essentially we just require that the bump moves modify in a transitive way  the
admissible boundary configurations $\fG'$. Then we can make the following informal conjecture

\begin{conj}\label{simple-proof} There is a substantially simpler proof of
  Proposition \ref{move-col}.
\end{conj}

Note that by using Proposition \ref{bump} (b), we can  generate a 2D equivalence
relation on the boundary configurations $\fG'$; so a preliminary task would be to get a `simple'
proof that these configurations are 2D equivalent to each other.

Assuming a (as much as possible) satisfactory solution of the conjecture, it is clear how to prove
Theorem \ref{bI}. Given two branched ideal triangulations $(T_1,b_1)$ and $(T_2,b_2)$ of $\hat M$
with boundary configurations $\fG_1$ and $\fG_2$, we perform on both a positive branched triangular
$0\to 2$ move (i.e. a bubble move in terms of dual spines) to get $(T'_1,b'_1)$ and $(T'_2,b'_2)$
representing points in $\Ss^{id}(M',\fG'_1)$ and $\Ss^{id}(M',\fG'_2)$ respectively.
By Proposition \ref{move-col} there is a composite ideal $b$-transit $(T'_1,b'_1)\Rightarrow (\tilde T_1, \tilde b_1)$
where this last represents a point in $\Ss^{id}(M', \fG'_2)$. Then by using Proposition \ref{pont} we get a composite
ideal $b$-transit  $(\tilde T_1, \tilde b_1)\Rightarrow (T_2',b_2')$. This is already an alternative proof of
Theorem \ref{bC}. Now to get Theorem \ref{bI} we have just to apply the arch and related constructions to undo the bubble moves.

Finally, given any naked ideal triangulation $T$ of $\hat M$ (as we know it happens that it does not carry any branching),
by theorem 3.4.9 of \cite{BP2} there is a composite naked transit $T \Rightarrow T'$ (entirely composed by positive $2\to 3$
moves) such that $T'$ carries a branching. So a proof of the naked ideal connectivity results would be derived from the branched one. 

\section{On $\Nn\Aa\Bb^{id}(*)$}\label{NA}
In this section we assume that the $3$-manifold $M$ is oriented.
The non ambiguous transit equivalence of pre-branched  ideal triangulations $(T, \omega)$ of $\hat M$
has been widely studied in \cite{NA}. An important remark is that every pre-branching determines
a branching  $\partial \omega$ on $\partial \tilde T$ (with the notations of Section \ref{verticalF}). 

\begin{figure}[ht]
\begin{center}
 \includegraphics[width=10cm]{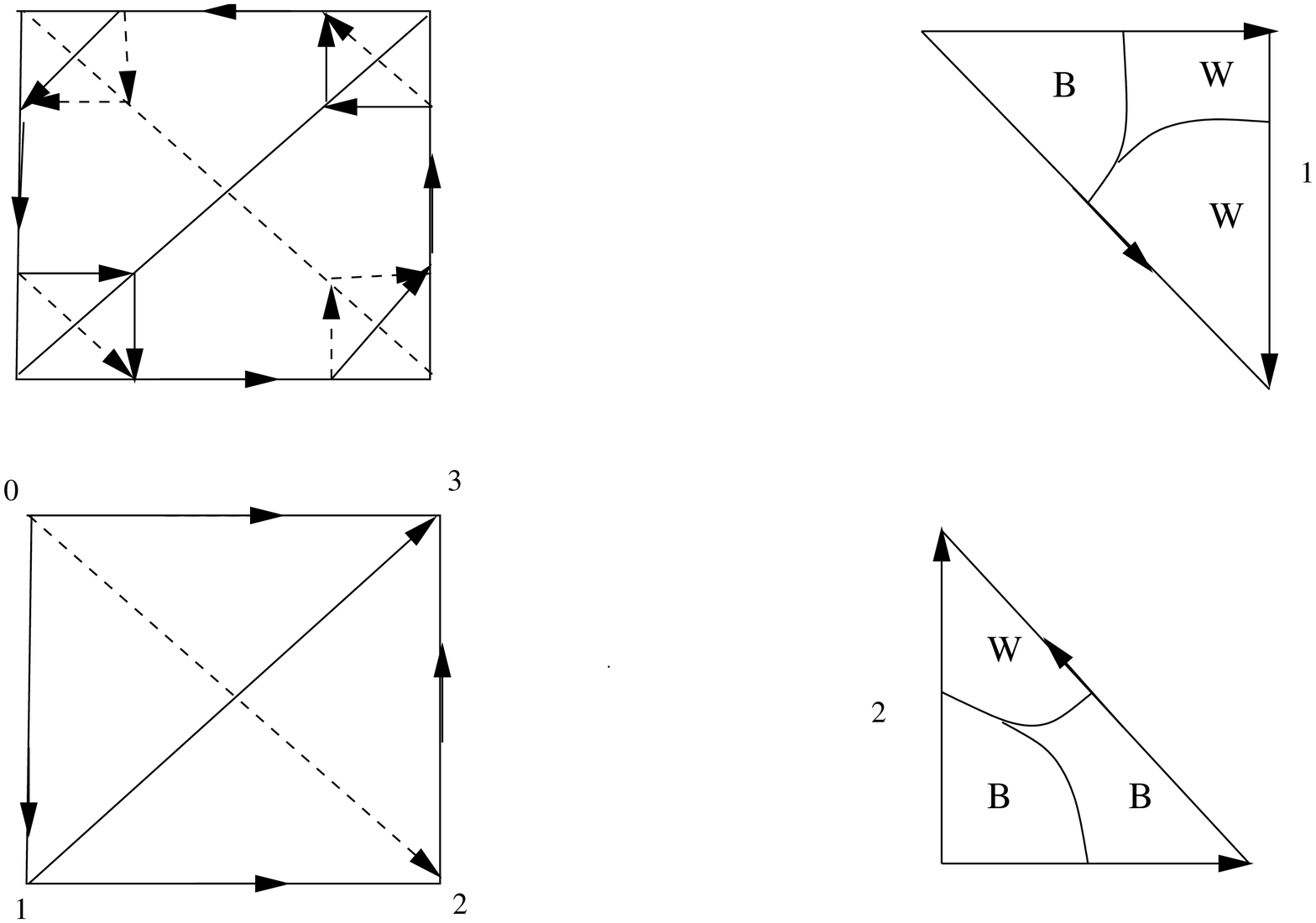}
\caption{\label{Vpuzzle} Tiles for $(\Vv,\partial \Vv)$.}
\end{center}
\end{figure}

This is illustrated in Figure \ref{Vpuzzle}. The picture shows:
 \begin{enumerate}
 \item A pre-branched tetrahedron $(\Delta,\omega)$ and a branched $(\Delta,b)$
 such that $\omega=\omega_b$.
 \item The four branched triangles in $(\partial \bar \Delta, \omega)$ (we stress again 
 that they only depend on the pre-branching $\omega$).
 \item The bicoloring of $t_1$ and $t_2$ determined by  $\Vv(\bar \tau,b)$.
 By specifying the branching of $t_j$ by labelling the vertices by $\{0,1,2\}$
 as usual, we see the black portion of $t_1$ contains only $v_0$, while the
 one of $t_2$ contains both $v_0$ and $v_1$; $v_2$ is always in the white
 portion. This qualitative behavior does not depend on the choice of  $(\Delta,b)$
 such that $\omega = \omega_b$, in particular on the sign $*_b$ (equal to $1$ in the picture).
 \end{enumerate}

\smallskip

The  2D version of the sliding equivalence  of branched triangulations of  
surfaces can be developed as well. A main invariant of 3D non ambiguous
classes of pre-branchings is the sliding equivalence class of $(\partial \tilde T, \partial \omega)$, hence
ultimately the associated pair of  oriented vertical and horizontal transvese foliations (with isolated singularities) 
on $\partial \tilde T$.  In fact we are considering here the branched triangulations of surfaces all together, with arbitrary
number of vertices, under both sliding diagonal exchanges and $1\to 3$ stellar moves (see \cite{NA} and \cite{2D}). 

If $(T,b)$ is branched, we can apply the above constructions to the induced pre-branching $\omega_b$.
Firstly we note that  the curve $X_\Vv$ of tangency lines of the vertical foliation $\Vv=\Vv_{T,b}$ is smoothly embedded 
into the oriented train track $\theta_b$ in $\partial M$ dual to $(\partial \tilde T, \partial \omega_b)$.
The vertical foliations  on $(\partial \tilde T, \partial \omega_b)$
can be recovered by means of the horizontal foliation $\Hh=\Hh_{T,b}$ of $(M,\partial M)$.
We outline this construction. 
\smallskip

{\bf (The ``maws'')} Consider $(T,b)$, $(\Sigma,b)$,
$(\Sigma^*,b^*)$ as in Section \ref{verticalF}. 

(1) We define a branched singular foliation $\mu$ on
$(\Sigma^*,b^*)$ (hence on $(\Sigma,b)$).  First we define $\mu$ along
a branched regular neighbourhood $\Nn$ of Sing$(\Sigma^*)$ in
$\Sigma^*$, in such a way that:
\begin{itemize}
\item It is traversing $\Nn$;
\item It is transverse to Sing$(\Sigma^*)$ and points everywhere
  toward the ``maw'' i.e. towards the $2$-regions of $(\Sigma^*,b^*)$
  whose orientation induces the {\it non} prevailing orientation on
  Sing$(\Sigma^*)$;
\item It has simple tangency points at $\partial \Nn$. 
\end{itemize}
\medskip

The closure of every component of $\Sigma^* \setminus \Nn$ is a
$2$-disk $R$. For every $\partial R$ there is an even number $2t(R)$
of tangency points of the partial foliation $\mu$ already
constructed. We define the index
$$d(R) = 1 - t(R) \ . $$  
Then we can extend $\mu$ to a singular foliation defined on the whole
of $\Sigma^*$. The singular set is contained in the set of
intersection points $x_R$ of the disks $R$ with the dual edges in $T$.
The smooth local model of $\mu$ at every such a point $x_R$ is either
like the vertical foliation at $0$ of the quadratic differential
$z^{-2d(R))}dz^2$, or is given by the gradient of $\pm(x^2 +
y^2)$. Hence $x_R$ is singular for $\mu$ if and only if $d(R)\neq
0$. Let us call $\mu$ the {\it maw foliation} of $\Sigma^*$ (hence on
$\Sigma$); it is uniquely determined up to homotopy through foliations
having the same properties along Sing$(\Sigma^*)$ and at the points
$x_R$'s. 
Notice that
  $$\chi(M) =  \frac{1}{2}\chi(\partial M) =  
\sum_{R\in \Sigma^{(2)}} d(R) \ $$
moreover, the cellular cochain $E(\Sigma,b) \in C^2(\Sigma; \Z)$ which
assigns the value $d(R)$ to every $(R,b)$ as above, represents the
Euler class of the oriented $2$-plane distributions associated to
$\Hh$.

\bigskip
 
(2) Take a normal embedding $M\subset F(T,b)$ We can consider now the
singular foliation $\hat \mu$ of $F = F(T,b)$ whose leaves are of the
form $l \times {t}$, $l$ being a leaf of $\mu$.  Hence for every
singular point $x$ of $\mu$, there is the vertical singular segment
$\{x\} \times I$ of $\hat \mu$.  Let $(\Vv, \Hh)$ be a couple of
vertical/horizontal foliations of $M$ constructed sofar.  For every
singular point $x$ of $\mu$, the leaf $\{x\} \times I \subset
M^\circ$, and has endpoints which belong to boundary leaves of both
$\Hh^*$ and $\Hh$.  Then $\hat \mu$ ``restricts'' to every component
of $\partial M \setminus \Uu_\Vv$.  This extends to a singular
foliation $\partial_s \mu$ (called the {\it singular boundary maw
  foliation}) of the whole of $\partial M$ whose singular set consists
of the union of the isolated maw singular points with the tangency
line $X=X_\Vv$. Finally the {\it boundary maw} $\partial \mu$ is obtained by
inverting the $\partial_s \mu$ orientations on the white components of
the bicoloring $\partial \Vv$ of $\partial M$.  It turns out that the
orientation conflict disappears at every component $C$ of $X_\Vv$ and
eventually $\partial \mu$ is traversing every annulus $\Uu_C$ without
boundary tangency points, and only the isolated maw singular points
survive; $\partial \mu$ is uniquely determined up to homotopy
through foliations positively transverse to $X_\Vv$ and which are
isotopic at the singular points. 
Finally we can state
\begin{prop}\label{2Dmu=v}   The singular foliation  $\partial \mu = \partial \mu_{T,b}$ coincides
with the vertical foliation carried by $(\partial \tilde T, \partial \omega_b)$.
\end{prop}

Essentially this holds `by construction'; we omit the details of the verification.
Now we can characterize the elementary non ambiguous ideal $b$-transits in terms of 
a preserved boundary configuration; the invariance of the maw besides the
boundary bicoloring characterizes the non ambiguous within the whole set of sliding transits.

\begin{prop}\label{non amb} Let $(T,b)\to (T',b')$ be an ideal ($2\leftrightarrow 3$ or $0\leftrightarrow 2$)
$b$-transit. Then the following facts are equivalent to each other:
\medskip

(1) It is non ambiguous;
\medskip

(2) $(\partial \Vv, \partial \mu)_{(T,b)}= (\partial \Vv, \partial \mu)_{(T',b')}$
up to isotopy. 
\smallskip

Moreover,  $(\partial \Vv, \partial \mu)_{(T,b)}$ is invariant under the non ambiguous
ideal branched equivalence.
\end{prop}

This can be checked case by case. 

We are going to finish by pointing out a further branched $na$ invariant
based on the horizontal foliation $\Hh=\Hh_{T,b}$.
We denote by $H^C_2(\Sigma; \R)$ the cellular singular homology,
provided that every region of the spine $\Sigma$ is oriented by the branching $b$.
Then $H^C_2(\Sigma; \R)$ is isomorphic to the singular homology $H_2(M; \R)$
and coincides with the space of $2$-cycles $Z^C_2(\Sigma; \R)$. 
Every $z\in Z^C_2(\Sigma; \R)$ consists in giving each $b$-oriented region
$R$ of $\Sigma$ a weight $z(R)\in \R$ in such a way that  the three weights around 
every edge of Sing$(\Sigma)$ verify
a switching condition of the form $z(e_0)=z(e_1)+z(e_2)$.
These cycles transit along every ideal $b$-transit, so that for every composite
$b$-transit $(T,b)\Rightarrow (T',b')$ it is defined an isomorphism
$$\alpha: Z^C_2(\Sigma;\R)\to Z^C_2(\Sigma';\R) \ . $$
Set
$$\Mm=\Mm_{(T,b)} = \{z\in Z^C_2(\Sigma; \R)| \ \forall 
R , \ z(R)\geq 0\}; \
\Mm^+ = \{z\in Z^C_2(\Sigma;\R)| \ \forall 
R, \ z(R)> 0\}\  . $$
Every $z\in \Mm$ can be interpreted as a {\it transverse measure}
on the horizontal foliation $\Hh$.
By taking into account the arbitrary choices in the realizations of 
$\Hh)$ we radily have

\begin{prop}\label{same measure}
  (1) For every $z\in \Mm$, the measured foliation $(\Hh,z)$
  is uniquely detemined up to {\rm measure equivalence}. 
    
  (2) If we denote by $\Mm(\Hh)$ the set of transverse measures 
on $\Hh$ up to measure equivalence, then the above correspondence 
well defines a map
$$ \mG=\mG_{(T,b)}: \Mm_{(T,b)}  \to \Mm(\Hh)  \  . $$
\end{prop}

\cvd

After a look at the $na$-transits we readily have

\begin{prop}
  If $(T,b)$ and $(T',b')$ are $na$-transit equivalent, then the maps
  $\mG_{(T,b)}$ and $\mG_{(T',b')}$ have the same image. More
  precisely, there is a bijection $\alpha: \Mm_{(T,b)}\to
  \Mm_{(T',b')}$ such that $\alpha(\Mm^+_{(T,b)})= \Mm^+_{(T',b')}$
  and $\mG_{(T,b)}= \mG_{(T',b')}\circ \alpha$.
\end{prop}

\cvd

\end{document}